\documentclass [12pt,a4paper]{article}

\usepackage{amsmath}
\usepackage{amsthm}
\usepackage{amssymb}

\usepackage{amscd}
\usepackage{amsfonts}
\usepackage{amsbsy}

\usepackage{graphicx}
\usepackage{indentfirst, latexsym, bm,amssymb}
\usepackage{bbding}

\advance\textwidth by +1.0in \advance\textheight by +1.0in
\advance\oddsidemargin by -0.5in \advance\evensidemargin by -1.0in
\advance\topmargin by -0.5in
\bibliography{bibfile}

\newtheorem {theorem} {Theorem}

\newtheorem {corollary} [theorem]{Corollary}
\newtheorem {lemma}  [theorem]{Lemma}

\def\v{\varepsilon}

\title{\large\bf Analytic Integrable
Systems: Analytic Normalization and Embedding Flows}
\author{\normalsize\bf\sc Xiang Zhang\footnote{\small
The author is partially supported by the NNSF of China grants 10831003 and 11271252, by RFDP of Higher Education of China grant 20110073110054, and by FP7-PEOPLE-2012-IRSES-316338 of Europe.}\\
\normalsize\it Department of Mathematics, and MOE-LSC, Shanghai Jiao Tong
               University,  \\ \normalsize\it Shanghai 200240,
               People's Republic of China.
      \\ \normalsize  E-mail: xzhang@sjtu.edu.cn\\}
\date{}
\begin{document}
\maketitle
\begin{abstract}
\noindent In this paper we mainly study the existence of analytic normalization and the normal form of finite dimensional complete analytic integrable dynamical systems. More details, we will prove that any
 complete analytic integrable diffeomorphism $F(x)=Bx+f(x)$ in $(\mathbb C^n,0)$ with $B$ having
eigenvalues not modulus $1$ and
$f(x)=O(|x|^2)$ is locally analytically conjugate to its normal form. Meanwhile, we also prove that
any  complete analytic integrable differential system $\dot x=Ax+f(x)$
in $(\mathbb C^n,0)$ with $A$ having nonzero eigenvalues and
$f(x)=O(|x|^2)$ is locally analytically conjugate to its normal form.
Furthermore we will prove that any complete analytic integrable
diffeomorphism defined on an analytic manifold can be embedded in
a complete analytic integrable flow. We note that parts of our results are the improvement of
Moser's one in {\it Comm. Pure Appl. Math.} 9$($1956$)$, 673--692 and of Poincar\'e's one
in {\it Rendiconti del circolo matematico di Palermo} 5$($1897$)$, 193--239.
These results also improve the ones in {\it J. Diff. Eqns.} 244$($2008$)$, 1080--1092 in the sense that
the linear part of the systems can be nonhyperbolic, and the one in {\it Math. Res. Lett.} 9$($2002$)$, 217--228
in the way that our paper presents the concrete expression of the normal form in a restricted case.

\hskip0.00011mm


\hskip0.00011mm

\end{abstract}

\bigskip
\section*{ \normalsize 1. Introduction and statement of the main results}

\setcounter{section}{1} \setcounter{equation}{0}

\noindent The study on the existence of analytic normalization for
an analytic dynamical system to its normal form has a long
history, which can be traced back to Poincar\'e and even earlier
(see e.g. \cite{Po,Mo,Br,XZ08,St09}). For analytic dynamical
systems, if their analytically equivalent normal forms are known,
it will be useful to study the dynamics of the original systems.
It is well--known, see e.g. \cite{Po,It,It1,Zu,XZ08,It09} and the
references therein, that the existence of analytic normalizations
for analytic vector fields to their normal forms is also strongly
related to the existence of analytic first integrals of
analytic vector fields.

The aim of this paper is to settle the problems on the existence
of analytic normalizations for analytic integrable diffeomorphisms to their normal forms and also for analytic integrable
vector fields to their normal forms.

For a diffeomorphism $F(x)$ defined in $(\mathbb
C^n,0)$, a function $V(x)$ is an {\it analytic first integral} of
$F(x)$ if it is analytic and satisfies $V(F(x))=V(x)$ for all
$x\in (\mathbb C^n,0)$. The diffeomorphism $F(x)$ is {\it
analytic integrable} in $(\mathbb C^n,0)$ if it has $n-1$
functionally independent analytic first integrals.
We should mention that the notion of integrable diffeomorphisms appeared only in recent years, see for instance \cite{CGM1,CGM,CN}.
As we know, there is no a notion of integrability on diffeomorphisms defined in the broad sense as that extended by Bogoyavlenski \cite{Bo} in 1998
for vector fields.

Denote by $M_n(\mathbb C)$ the set of square matrices of order $n$ with entries in $\mathbb C$. Let $B\in M_n(\mathbb C)$ and $\mu=(\mu_1,\ldots,\mu_n)$ be the eigenvalues of $B$.
Recall that
\begin{itemize}

\item A diffeomorphism or a formal series $F(x)=Bx+f(x)$ is in {\it normal
form} if $B$ is in the Jordan normal form and the
nonlinear term $f(x)$ consists of only resonant monomials. A
monomial $x^me_j$ in the $j$th component of $f(x)$ is {\it
resonant} if $\mu^m=\mu_j$, where $e_j$ is the unit vector with its $j$th component
equal to 1 and others vanishing.

\item A diffeomorphism (or a formal series) $G(y)$ is a {\it normal form } (or a {\it formal normal form}) of a diffeomorphism $F(x)$ if
$G(y)$ is in normal form and $F(x)$ and $G(y)$ are conjugate, i.e., there is a transformation tangent to identity $y=\Phi(x)=x+\phi(x)$
with $\phi(x)$ containing only higher order terms such that $G\circ \Phi(x)=\Phi \circ F(x)$. The conjugacy $y=\Phi(x)$ is called a {\it normalization} from $F(x)$ to $G(y)$. Furthermore

\begin{itemize}

\item  $y=\Phi(x)$ is an {\it analytic normalization} of $F(x)$ if $\Phi(x)$ is analytic. In this case we call $G(y)$ an {\it analytically equivalent normal form} of $F(x)$.

\item  $y=\Phi(x)$ is a {\it distinguished normalization} of $F(x)$ if $\phi(x)$ contains only nonresonant term, i.e. its monomial $x^me_j$ in the $j$th component of $\phi(x)$ are  {\it nonresonant} in the sense that $\mu^m\ne 1$.

\end{itemize}
\end{itemize}
We remind readers the difference between the resonances of diffeomorphisms and the transformations.

Our first main result of this paper provides more information on
analytic integrable diffeomorphisms than the existence of
analytic normalization. Before stating the results, we introduce a
notation. Let $\mu=(\mu_1,\ldots,\mu_n)$ be the eigenvalues of
$B$. Set
\[
\mathcal D=\{m\in\mathbb Z_+^n; \,\, \mu^m=1,|m|\ge 2\},
\]
namely {\it resonant set} of $B$, and denote by $d_\mu$ the rank
of the resonant set. The elements of $\mathcal D$ are also called
{\it resonant lattices}. An element $m\in\mathcal D$ is {\it
simple} if it cannot be divided by a positive integer no less than
$2$.

\begin{theorem}\label{t2}
For a diffeomorphism $F(x)=Bx+f(x)$ defined in $(\mathbb C^n,0)$  a
neighborhood of $0$ in $\mathbb C^n$ with $f(x)=O(|x|^2)$ and $B$ having at least one eigenvalue not on
the unit circle of $\mathbb C$, then $F(x)$ is analytic
integrable if and only if the following statements hold.
\begin{itemize}
\item[$(a)$] the resonant set has the rank $d_\mu=n-1$.

\item[$(b)$] $F(x)$ is conjugate to its normal form
of type
\[
G(y)=(\mu_1 y_1(1+p_1(y)),\ldots,\mu_ny_n(1+p_n(y)),
\]
by a distinguished analytic normalization, where
$\mu=(\mu_1,\ldots,\mu_n)$ is the $n$--tuple of eigenvalues of
$B$, and $p_1(y),\ldots,p_n(y)$ are analytic and  satisfy the
equations
\[
(1+p_1(y))^{m_{k1}}\ldots(1+p_n(y))^{m_{kn}}=1, \quad\mbox{ for }
k=1,\ldots,n-1,
\]
with $m_k=(m_{k1},\ldots,m_{kn})\in\mathcal D$, $k=1,\ldots,n-1$,
being $n-1$ linearly independent simple resonant lattices.

\end{itemize}
\end{theorem}

We remark that this last result is a correction and improvement of the one given in \cite{XZ08}. From this last theorem and the following Lemma \ref{l24}(c) we can get easily the following.
\begin{corollary}\label{c1}
For analytic integrable diffeomorphism $F(x)=Bx+f(x)$ with
$f(x)$ nonlinear, if the orbits of the normal form system of
$F(x)$ and of the linear one $\mu y=(\mu_1y_1,\ldots,\mu_ny_n)$
start at the same generic point, then the full orbits will be contained in
the same orbit of an analytic vector field, where
$\mu=(\mu_1,\ldots,\mu_n)$ is the $n$--tuple of eigenvalues of $B$
and not all of $\mu_i$'s on the unit circle of $\mathbb C$.
\end{corollary}

Recall that a {\it generic point} is the one which is located in a full Lebesgue measure subset of $(\mathbb C^n,0)$.

As we know, for higher dimensional local analytic diffeomorphisms
the existence of analytic normalization is solved only  for the
diffeomorphisms having their linear parts with eigenvalues
$\mu=(\mu_1,\ldots,\mu_n)$ either all larger (or smaller) than $1$
in modulus, or satisfying $|\mu^m-\mu_s|\ge c|m|^{-\nu}$ for all
$s=1,\ldots,n$, $|m|\ge 2$, with $c,\nu>0$ given constants. The
former result is called {\it Poincar\'e--Dulac theorem}, and the
latter is called {\it Siegel theorem}, see e.g. $\S 25$ of
\cite{Ar}. We note that in the Siegel theorem the eigenvalues are
nonresonant. In our case the eigenvalues of the linear parts of the diffeomorphisms can be resonant,
and their modulus can have part of them  larger than $1$ and have also other part of them less than $1$.

We now turn to the study of the problem on the existence of
analytic normalizations for  analytic integrable vector fields to
their normal forms. We will see
that this problem is simpler than that for analytic integrable diffeomorphisms.

Consider the analytic
differential system
\begin{equation}\label{e1}
\dot x=Ax+f(x),\qquad x\in \left(\mathbb C^n, 0\right),
\end{equation}
where $A\in M_n(\mathbb C)$, and $f(x)=O(|x|^2)$
is a  vector-valued analytic function in $(\mathbb C^n,0)$. We say that system
\eqref{e1} is {\it locally complete analytic integrable} in $(\mathbb
C^n,0)$ if it has $n-1$ functionally independent analytic first
integrals in $(\mathbb C^n,0)$. An {\it analytic first integral}
of system \eqref{e1} is a nonconstant analytic function $H(x)$
defined in $(\mathbb C^n,0)$ satisfying $\langle\nabla
H(x),Ax+f(x)\rangle=0$ in $(\mathbb C^n,0)$, where $\langle
\cdot,\cdot\rangle$ denotes the inner product of two vectors in
$\mathbb C^n$, and $\nabla $ represents the gradient of a function
with respect to $x$. The $n-1$ analytic first integrals are {\it
functionally independent} if their gradients as vectors in
$\mathbb C^n$ are linearly independent on an open dense subset of
$(\mathbb C^n,0)$.

Let $\lambda=(\lambda_1,\ldots,\lambda_n)$ be the $n$--tuple of
eigenvalues of the matrix $A$. Set
\[
\mathcal R_\lambda:=\left\{m=(m_1,\ldots,m_n)\in\mathbb Z_+^n;
\,\,\,\langle m,\,\lambda\rangle=0,\,\,|m|=m_1+\ldots+m_n\ge
2\right\},
\]
where $\mathbb Z_+$ denotes the set of nonnegative integers. We
call $\mathcal R_\lambda$  {\it resonant set} of $\lambda$, and
its elements  {\it resonant lattices}. Denote by $r_\lambda$ the
rank of vectors in the set $\mathcal R_\lambda$. Obviously, if
$\lambda\ne 0$ then $r_\lambda\le n-1$. In what follows we assume
without loss of generality that $A$ is in the lower triangular Jordan
normal form.

In \cite{XZ08} we have proved the following result.  {\it Assume
that the origin of system \eqref{e1} is nondegenerate, i.e. no
eigenvalues equal to zero, and that the matrix $A$ is
diagonalizable. Then system \eqref{e1} has $n-1$ locally
functionally independent analytic first integrals if and only if
$r_\lambda=n-1$, and system \eqref{e1} is analytically equivalent
to its distinguished normal form $ \dot y_i=\lambda_iy_i(1+g(y))$,
$i=1,\ldots,n$, by an analytic normalization, where $g(y)$,
without constant term, is an analytic function of $y^m$ with
$m\in\mathcal R_\lambda$ and $(m_1,\ldots,m_n)=1$, i.e. they have
no common factor.} Recall by definition  that $y^m=y_1^{m_1}\ldots
y_n^{m_n}$ for $y=(y_1,\ldots,y_n)$ and $m=(m_1,\ldots,m_n)$.

We say that system \eqref{e1} is in {\it normal form} if $A$ is in
Jordan normal form and the Taylor series of $f(x)$ consists of
only resonant monomials. A monomial $x^me_j$ modulo coefficient in
the $j$th component of $f(x)$ is {\it resonant} if
$\lambda_j=\langle m,\lambda\rangle$, where $m\in\mathbb Z_+^n$,
$|m|\ge 2$.

A system
\begin{equation}\label{e1.nf1}
\dot y=Ay+g(y),
\end{equation}
with $g(y)$ containing only higher order terms is a {\it normal form} of
\eqref{e1} if system \eqref{e1.nf1} is in normal form, and there is a change of variables tangent to identity $y=\Psi(x)$ transforming
system \eqref{e1} to \eqref{e1.nf1}. The transformation $y=\Psi(x)$ is called a {\it normalization}. If the transformation contains only nonresonant terms, then it is called {\it distinguished normalization}.
Correspondingly, the normal form is called a {\it distinguished
normal form}. Recall that a monomial $x^me_j$ in the transformation  is {\it nonresonant} if $\langle
m,\lambda\rangle\ne 0$.

Related to the above results we posed in Remark 2 of \cite{XZ08}
the following open problem: {\it if system \eqref{e1} has the
origin as a degenerate singularity, and has $n-1$ locally
functionally independent analytic first integrals in a
neighborhood of the origin, is system \eqref{e1} locally
analytically equivalent to its distinguished normal forms?}

In this paper we will give a positive answer to this problem. The following is our second main result of this paper.

\begin{theorem}\label{t1}
Assume that $n\ge 2$ and $\lambda\ne 0$, i.e. $A$ has at least one
eigenvalues not equal to zero. Then system \eqref{e1} has $n-1$
functionally independent analytic first integrals in $(\mathbb
C^n,0)$ if and only if the rank of the resonant set is equal to
$1$, i.e. $r_\lambda=n-1$, and system \eqref{e1} is analytically
equivalent to its distinguished normal form
\begin{equation}\label{e2}
\dot y_i=\lambda_iy_i(1+g(y)),\qquad i=1,\ldots,n,
\end{equation}
by an analytic normalization tangent to the identity, where
$g(y)$, without constant term, is an analytic function of $y^m$
with $m\in\mathcal R_\lambda$ and $(m_1,\ldots,m_n)=1$.
\end{theorem}

This last theorem characterizes the  complete  analytic integrable
differential systems. It provides not only the existence of
analytic normalization for analytic  integrable differential
systems but also the concrete expression of their normal forms.
Theorem \ref{t1} also implies that the linear part of an
analytic integrable differential system must be diagonalizable
provided that it has at least one nonzero eigenvalue.
We can improve the corresponding result given in \cite{XZ08}, because
we find that in the degenerate case there also does not appear small divisors
in the normalization.

We should mention that one dimensional differential equation, if
it is nontrivial, has no nonconstant first integrals. For
higher dimensional systems, if all eigenvalues $\lambda$ of $A$
are equal to zero, Theorem \ref{t1} cannot be applied, see the
following examples.

\noindent{\it Example }1. The planar system
\[
\dot x= y+f(x,y),\quad \dot y=0,
\]
with $f(x,y)$ an analytic function without linear and constant
terms, has the analytic first integral $H(x,y)=y$. Also the planar
analytic system
\[
\dot x=y^{p+1},\quad \dot y=-x^{q+1}, \qquad p,q\in \mathbb N,
\]
has the analytic first integral
$H(x,y)=x^{q+2}/(q+2)+y^{p+2}/(p+2)$.  But these two analytic
integrable systems cannot be transformed to systems of form \eqref{e2}
by an invertible analytic change of coordinates, because the
linear part of these two systems have only zero eigenvalues. Also
we note that according to the definition of resonance the
nonlinear parts of these last two systems are all resonant.

These last examples show that the condition $\lambda\ne 0$  is
necessary for Theorem \ref{t1}. So we have the following

\noindent {\it Open problem} 1. Assume that the linear part of
system \eqref{e1} has all eigenvalues vanishing.
\begin{itemize}
\item{} What is the normal form that an analytic integrable
system \eqref{e1} can have? \item{}  Is an analytic integrable
system \eqref{e1} analytically equivalent to its normal form?
\end{itemize}

Comparing with Theorem \ref{t1}, we must mention the work of Zung \cite{Zu1}, in which author studied the existence of analytic normalization for obtaining the Poincar\'e--Dulac
normal form of an analytic integrable differential systems in the broad sense. In 1998
Bogoyavlenski \cite{Bo} extended the classical complete integrable differential systems including the Liouvillian integrable Hamiltonian systems as follows: a local analytic (or smooth) vector field in $(\mathbb C^n,0)$ is {\it analytic } (or {\it smooth}) {\it integrable in the broad sense} if for some natural number $q$ ($1\le q\le n$), there exist $q$ locally analytic (or smooth) vector fields $X_1=X$, $X_2,\ldots,X_q$ and $n-q$ locally analytic (or smooth) functions $f_1,\ldots,f_{n-q}$ in $(\mathbb C^n,0)$ such that
\begin{itemize}
\item[$i)$] the $q$ vector fields $X_1,\ldots,X_q$ commute pairwise and are linearly independent almost everywhere;
\item[$ii)$] $f_1,\ldots,f_{n-q}$ are common first integrals of $X_1,\ldots,X_q$ and are functionally independent almost everywhere.
\end{itemize}
Using a geometric method Zung \cite{Zu1} in 2002 proved that
any analytic integrable differential system in the broad sense is analytically conjugate to its normal form.
From the existence of analytic normalization point of view, Zung's result  contains our Theorem \ref{t1} as a special case. But our result can present the concrete expression of the normal form of the system. Since we do not study this broad sense integrability, in the following when we say integrable systems, we always mean the complete integrable systems which were defined above.

We note that  the distinguished normal form of an analytic integrable vector field has a concrete expression,
but we cannot present the exact expression of their normal forms  for analytic integrable diffeomorphisms. In fact, they depend on
the resonant lattices. This can be seen from the following concrete examples, which illustrate some applications of
Theorem \ref{t2}.

\noindent{\it Example} 2. According to Theorem \ref{t2}, the two
dimensional analytic integrable diffeomorphism $(\frac 12
x+g_1(x,y),\,2y+g_2(x,y))$ is locally analytically conjugate to
its distinguished normal form
\[
F_1(x,y)=\left(\frac 12
x\left(1-\frac{\varphi(xy)}{1+\varphi(xy)}\right),\,\, 2
y(1+\varphi(xy))\right),
\]
with $\varphi(z)$ an analytic function satisfying $\varphi(0)=0$.
$F_1(x,y)$ has  the same first integral $H(x,y)=xy$ as the linear
diffeomorphism $L_1(x,y)=(\frac 12 x,2y)$.  But they cannot be
parallel as in the case of two dimensional vector fields.

The three dimensional analytic integrable diffeomorphism
\[
\left(e^{-5}x+g_1(x,y,z),\,\,e^2y+g_2(x,y,z),\,\,e
z+g_3(x,y,z)\right),
\]
is locally analytically conjugate to its
normal form
\[
F_2(x,y,z)=\left(e^{-5}x\frac{1}{(1+\psi)^{5/2}},\,\,e^2y(1+\psi),\,\,e
z(1+\psi)^{1/2}\right),
\]
which together with the linear diffeomorphism
$L_2(x,y,z)=(e^{-5}x,\,e^2 y,\,e z)$ have the functionally
independent analytic first integrals $H_1(x,y,z)=x^2 y^5$ and
$H_2(x,y,z)=xy^2z$, where
$\psi=\psi(w_1,w_2,w_3,w_4)=\psi(xy^2z,xyz^3,xz^5,x^2y^5)$ is an
analytic function in its variables and $\psi(0,0,0,0)=0$. We note
that $xyz^3$ and $xz^5$ are also analytic first integrals of $F_2$
and of $L_2$, and  they functionally depend on $xy^2z$ and
$x^2y^5$. But they cannot be represented in analytic functions of
$xy^2z$ and $x^2y^5$.

In addition, any three dimensional locally analytic integrable
diffeomorphism
\[
G(x,y,z)=(x+g_1(x,y,z)),\,y+g_2(x,y,z)),\,2z +g_3(x,y,z)),
\]
with $g_1,g_2,g_3$ nonlinear, is analytically conjugate to
$(x,\,y,\,2z(1+h(x,y)))$ with $h$ an analytic function in $x$ and
$y$.

The last three examples show that the normal forms of higher
dimensional analytic integral diffeomorphisms have more
complicated expressions than those of vector fields.

Now we briefly review the results on the existence of analytic
normalizations for analytic integrable vector fields to their
normal forms and also on the existence of analytic integrable
vector fields.

The following result, known as Poincar\'e normal form theorem of a
nondegenerate center, goes back to Poincar\'e and Lyapunov (see
e.g. \cite{Po}): a planar analytic differential system has the
origin as a nondegenerate center if and only if it is
analytically equivalent (via probably complex transformation of
variables and time rescaling) to
\begin{equation}\label{1.1}
\dot x=x(1+q(xy)),\quad \dot y=-y(1+q(xy)),
\end{equation}
where $q(u)$ is an analytic function in $u$ starting from the
terms of degree no less than 1. We note that this result is a
special case of our Theorem \ref{t1}, because in the case of
nondegenerate center the unique linearly independent simple
resonant lattices is $(1,1)$ under the complex coordinates. The
Poincar\'e normal form theorem has a corollary as follows: a
planar analytic differential system has the origin as an
isochronous center if and only if it is analytically equivalent to
\[
\dot u=-\omega v,\quad \dot v=\omega u,
\]
where $\omega$ is a nonzero constant. Moser \cite{Mo} showed
that a planar real analytic Hamiltonian system having the origin
as a hyperbolic saddle can be reduced to system \eqref{1.1} by a
real analytic area-preserving transformation of variables. This
shows that Theorem \ref{t1} is not only a generalization of the
above Poincar\'e's and Moser's results to higher dimensional
systems in nondegenerate cases, but also a generalization of
their results to degenerate cases. For example, the planar
analytic differential system $\dot x=0,\,\,\dot y=\lambda
y+o(|x,y|^2)$ with $\lambda\ne 0$ is analytically equivalent to a
system of form $\dot x=0, \,\,\dot y=y(\lambda+O(x))$ by Theorem
\ref{t1}, because the system has a functionally independent
analytic first integral. For general planar analytic differential
systems, Llibre {\it et al.} \cite{CGGL,CGGL1,GGL} characterized
their locally analytic integrability around a singularity with the
aid of normal forms. Closely related to analytic integrability of
planar differential systems, the existence of inverse analytic
integrating factors provides much more information on the dynamics
of the system (see e.g. \cite{EP09,GGG09,GM09}).

Our study on the existence of analytic normalization and the concrete expressions of the normal forms for complete analytic
integrable systems in $\mathbb C^n$ is also motivated by the study
on a similar problem for Hamiltonian systems, that is, on the
existence of analytically symplectic normalizations which transform
analytic integrable symplectic Hamiltonian systems to their
Birkhoff normal form.  Ito \cite{It, It1} solved this problem
under the restriction that the eigenvalues of linear part of
systems are nonresonant and simple resonant, respectively. Zung
\cite{Zu} completely solved this problem, and proved that
analytic Liouvillian integrable symplectic Hamiltonian system
is analytically symplectically equivalent to its Birkhoff normal
form by developing a new geometric method based on the toric
characterization of Birkhoff normalization. Recently Ito
\cite{It09} presented a relation between superintegrability of
Hamiltonian systems and the existence of analytic Birkhoff
normalization.

Here we mainly concern the analytic normalization for analytic
integrable systems. For general differential systems including
Hamiltonian ones, there are extensive studies on the existence of
analytic normalizations, we refer readers to the papers
\cite{St09, St}, the books \cite{Li00,IY08} and the references
therein. On the generic nonexistence of analytic normalization
for analytic differential systems, we refer readers to Siegel
\cite{Si} and P\'{e}rez-Marco \cite{Pe1, Pe}.

Theorems \ref{t2} and \ref{t1} show that both analytic
integrable diffeomorphisms and analytic integrable
vector fields have analytic normalizations. It motivates to think
whether an analytic integrable diffeomorphism can be embedded
in an analytic integrable autonomous vector field.

The following result was given in \cite{XZ08} (for a similar one,
see Cima {\it et al} \cite{CGM}).  {\it Any analytic
integrable volume-preserving diffeomorphism defined on an analytic
manifold $\mathcal M$ can be embedded in an analytic flow on
$\mathcal M$.}

Here we will release the restriction on volume--preserving of the
diffeomorphisms and will give a global proof on the given manifold
where the integrable diffeomorphisms are defined. Our result is
the following

\begin{theorem}\label{t3} Let $\mathcal M$ be a real or complex $n$--dimensional analytic
manifold. Then any real or complex analytic integrable
diffeomorphism defined on $\mathcal M$ can be embedded in an
analytic flow on $\mathcal M$.
\end{theorem}

Using the same method as that in the proof of Theorem \ref{t3}, we can get easily the
following result.

\begin{corollary}
Let $\mathcal M$ be an $n$--dimensional $C^k$ smooth manifold
with $k\in\{\mathbb N\}\cup\{\infty\}$. Then any $C^k$ smoothness
integrable diffeomorphism defined on $\mathcal M$ can be embedded
in a $C^{k-1}$ smoothness flow on $\mathcal M$, where
$\infty-1=\infty$.
\end{corollary}

These last results have  solved the open problem given in Remark 5
of \cite{XZ08}. Recall that for an analytic or smooth manifold
$\mathcal M$, a diffeomorphism $F:\mathcal M\rightarrow \mathcal
M$ can be {\it embedded in a flow $\phi_t$} if $\phi_1=F$ on
$\mathcal M$. The vector field $\partial \phi_t/\partial t|_{t=0}$
is called {\it embedding vector field} of $F$. Generally, a time
dependent vector field $\mathcal X(t,x)$ is an {\it embedding
vector field} of the diffeomorphism $F(x)$ if the latter coincides
with the time 1 map of solutions of $\mathcal X(t,x)$.

In one dimensional case there are rich results on the embedding
flow problem (see e.g. \cite{BT03, Li00} and the references
therein). In higher dimensional cases, Arnold \cite{Ar} posed the following
result without a proof that if $A\in M_n(\mathbb R)$ has a real logarithm $B$, i.e. $A=e^B$, then any $C^\infty$ local diffeomorphism
$f(x)=Ax+O(|x|^2)$ can be embedded in a $C^\infty$ periodic vector
field $\dot x=Ax+g(t,x)$ in $(\mathbb
R^n,0)$, where $g(t+1,x)=g(t,x)$ and $g(t,x)=O(|x|^2)$, for a proof see \cite[Lemma 16]{LLZ}). Kuksin and P$\rm\ddot{o}$schel \cite{KP}
proved the existence of $C^\infty$ or analytic periodic
embedding Hamiltonian vector fields for a class of nearly integrable  $C^\infty$ or analytic symplectic diffeomorphisms
defined in $(\mathbb R^{2m},0)$.

The problem on the existence of embedding flows or embedding autonomous vector fields becomes more
difficult, as mentioned by Arnold in \cite[p.200]{Ar}. Palis \cite{Pa}
proved that the diffeomorphisms admitting embedding flows are rare
in the Baire sense. In \cite{LLZ}, we proved that for a $C^\infty$ diffeomorphisms $F(x)=Ax+f(x)$ defined in $(\mathbb R^n,0)$
with $A$ having a real logarithm $B$ and $f(x)=O(|x|^2)$, if $A$ has no eigenvalues on the unit circle of $\mathbb C$ and the
eigenvalues of $B$ are not weakly resonant, then $F(x)$ can be embedded in a $C^\infty$ autonomous vector field. This result was recently
extended to Banach spaces \cite{XZ09}. For analytic
diffeomorphisms, as we know, Theorem 1.3 of \cite{XZ08} and part
of Theorem 1.4 of \cite{XZ10} are the only results on the
existence of analytic embedding flows in higher dimensional
spaces. Recently Zhang \cite{XZ2011} provided a simple proof to the result of \cite{LLZ} mentioned above and presented
some examples showing that the weakly nonresonant conditions of the real logarithm $B$ of $A$ is necessary.

Theorem \ref{t3} has solved the embedding flow problem for analytic integrable diffeomorphisms. But for nonintegrable analytic
diffeomorphisms the problem is still open.

\noindent{\it Open problem} 2. To characterize all analytic
diffeomorphisms which admit analytic embedding flows.

The paper is organized as follows. We first prove Theorem \ref{t2}
in Section 2. Because the proof of Theorem \ref{t1} is similar to
that of Theorem \ref{t2} and is easier, which will be given in Section
3, where we mainly concern the difference with that of
Theorem \ref{t2}. The last section presents the proof of Theorem
\ref{t3}.

\section{Proof of Theorem \ref{t2}}

\setcounter{section}{2}
\setcounter{equation}{0}\setcounter{theorem}{0}

\noindent  We separate the proof of Theorem \ref{t2} into
several lemmas. One of the main tools for proving the theorem is
the normal form theory. For doing so, we need an auxiliary result,
which will be used later on in different ways for several times.

\begin{lemma}\label{l20}
Let $\mathcal H_n^r(\mathbb C)$ be the linear space of
$n$--dimensional vector--valued homogeneous polynomials of degree
$r$ in $n$ variables with coefficients in $\mathbb C$. For $B,C\in
M_n(\mathbb C)$, we define a linear operator $\mathcal L_{B,C}$ on
$\mathcal H_n^r(\mathbb C)$ by
\[
(\mathcal L_{B,C}\phi)(x)=\phi(Bx)-C\phi(x),\qquad \phi\in\mathcal
H_n^r(\mathbb C).
\]
Then the spectrum, denoted by $\sigma(\mathcal L_{B,C})$, of
$\mathcal L_{B,C}$ is
\[
\sigma(\mathcal L_{B,C})=\{\mu^m-\kappa_j;\,\, m\in\mathbb Z_+^n,
|m|=r,j=1,\ldots,n\},
\]
where $\mu=(\mu_1,\ldots,\mu_n)$ and
$\kappa=(\kappa_1,\ldots,\kappa_n)$ are the $n$--tuples of
eigenvalues of $B$ and $C$, respectively.
\end{lemma}

\noindent{\it Proof}. The idea of the proof follows from that of Lemma 1.1 of
\cite{Bi} and of Lemma 4.5 of \cite{Li00}. Let $T,S\in M_n(\mathbb
C)$ be invertible. Set $\phi(x)=T\xi(x)$ and $x=Sy$. Then for
$\psi(y)=\xi(Sy)$ we have
\begin{eqnarray*}
(\mathcal
L_{B,C}\phi)(x)&=&T\xi(Bx)-CT\xi(x)=T\left(\xi(BSy)-T^{-1}CT
\xi(Sy)\right)\\
&=&T\left(\psi(S^{-1}BSy)-T^{-1}CT \psi(y)\right).
\end{eqnarray*}
Consider the linear operator
\begin{equation}\label{ebaffle}
(\mathcal L^*\psi)(y)=\psi(S^{-1}BSy)-T^{-1}CT \psi(y).
\end{equation}
Then the linear operators $\mathcal L_{B,C}$ and $\mathcal L^*$
have the same spectrum, because  $\psi(y)=T^{-1}\phi(Sy)$, and $S$
and $T$ are invertible. So without loss of generality we can
assume that the matrices $B$ and $C$ are in lower triangular
Jordan normal form.

\noindent{\it Case }1. $B$ and $C$ are diagonalizable. From \eqref{ebaffle},  we can assume without
loss of generality that $B$ and $C$  are diagonal. For any
monomial $h(x)=x^me_j\in \mathcal B=\{x^me_j;\,\, m\in\mathbb
Z_+^n, |m|=r,j=1,\ldots,n\}$ a base of $\mathcal H_n^r(\mathbb
C)$, where $e_j$ is the $j$th unit vector, we have
\[
\mathcal L_{B,C}(x^me_j)=(\mu
x)^me_j-\mbox{diag}(\kappa_1,\ldots,\kappa_n)x^me_j=(\mu^m-\kappa_j)x^me_j.
\]
This shows that the matrix expression of the linear operator
$\mathcal L_{B,C}$ under the base $\mathcal B$ is diagonal with
$\mu^m-\kappa_j$ being the elements  on the diagonal entries.
Hence $\mathcal L_{B,C}$ has the spectrum as stated in the lemma.

\noindent{\it Case} 2. At least one of $B$ and $C$ is not
diagonalizable. Choose $B(\varepsilon), C(\varepsilon)\in
M_n(\mathbb C)$ such that $B(\varepsilon)\rightarrow B$ and $
C(\varepsilon)\rightarrow C$ as $\varepsilon\rightarrow 0$ and
that $B(\varepsilon)$ and $C(\varepsilon)$ are both
diagonalizable. Let
$\mu(\varepsilon)=(\mu_1(\varepsilon),\ldots,\mu_n(\varepsilon))$
and
$\kappa(\varepsilon)=(\kappa_1(\varepsilon),\ldots,\kappa_n(\varepsilon))$
be the $n$--tuples of eigenvalues of the matrices $B(\varepsilon)$
and $C(\varepsilon)$, respectively. Then it follows from the proof
of Case 1 that the linear operator $\mathcal
L_{B(\varepsilon),C(\varepsilon)}$ has the spectrum $
\sigma(\mathcal
L_{B(\varepsilon),C(\varepsilon)})=\{\mu(\varepsilon)^m-\kappa_j(\varepsilon);\,\,
m\in\mathbb Z_+^n, |m|=r,j=1,\ldots,n\}$. Since the linear
operator $\mathcal L_{B,C}$ depends on $B$ and $C$ continuously,
and  $\mu(\varepsilon)\rightarrow \mu$ and
$\kappa(\varepsilon)\rightarrow \kappa$ as $\varepsilon\rightarrow
0$, we get the spectrum of $\mathcal L_{B,C}$ as stated in the
lemma. The proof is completed.

The first result is on the existence of formal normal form for
analytic diffeomorphims, which can be found in any book when it
introduces normal form theory (see e.g. \cite{Ar, Bi, Li00}).

\begin{lemma}\label{l21} The analytic diffeomorphism
$F(x)=Bx+f(x)$ is always formally conjugate to its distinguished
normal form by a formal transformation tangent to the identity.
\end{lemma}


\noindent{\it Proof}. We present its proof here because  which  will
be used in the proof of our other results. For the diffeomorphism
$F(x)=Bx+f(x)$, we can assume without loss of generality that $B$
is in lower triangle Jordan normal form. Because there always
exists an invertible linear conjugation $Cy$ such that $F(x)$ is
conjugate to $C^{-1} BCy+C^{-1}h(Cy)$.

Suppose that $F(x)=Bx+f(x)$ is conjugated to $G(y)=By+g(y)$ via a
formal conjugation tangent to the identity, i.e. it is of the form
$x=\Phi(y)=y+\phi(y)$ with $\phi(y)$ a formal series starting from
at least the second order term. Then we get from $F\circ
\Phi(y)=\Phi\circ G(y)$ that $g$ and $\phi$ satisfy
\begin{equation}\label{2.1}
\phi(By)-B\phi(y)=f(y+\phi(y))+\phi(By)-\phi(By+g(y))-g(y).
\end{equation}
Expanding $h\in\{f,g,\phi\}$ in Taylor series gives
\[
h(x)=\sum\limits_{i=l}\limits^{\infty}h_i(x),
\]
where $h_i$ is a vector-valued homogeneous polynomial of degree
$i$, and $l$ is the degree of the lowest order homogeneous
polynomial in the Taylor expansion of $f(x)$. Then we get from
\eqref{2.1} that
\begin{equation}\label{2.2}
\phi_s(By)-B\phi_s(y)=[f]_s+[\phi]_s-g_s(y),\qquad s=l,l+1,\ldots,
\end{equation}
where $[f]_s$ and $[\phi]_s$ are inductively known vector--valued
homogeneous polynomials in $y$ of degree $s$ obtained by
re--expanding $f(y+\phi(y))$ and $\phi(By)-\phi(By+g(y))$ in
Taylor series in $y$, respectively. In fact, we have $[\phi]_l=0$.

Recall that $\mathcal H_n^s(\mathbb C)$ is the linear space formed
by $n$--dimensional vector--valued complex homogeneous polynomials
of degree $r$ in $n$ variables. Define a linear operator $\mathcal
L_{B,s}: \mathcal H_n^s(\mathbb C)\rightarrow \mathcal
H_n^s(\mathbb C)$ by
\[
\mathcal L_{B,s}\phi(y)=\phi(By)-B\phi(y)\quad \mbox{ for }
\phi\in \mathcal H_n^s(\mathbb C).
\]
We get from Lemma \ref{l20} that the spectrum of $\mathcal
L_{B,s}$ is
\[
\left\{\prod\limits_{i=1}\limits^{n}\mu_i^{m_i}-\mu_j;\,
m_i\in{\mathbb Z}_+,\,\sum\limits_{i=1}\limits^{n}m_i=s,\,
j=1,\ldots,n\right \},
\]
where $(\mu_1,\ldots, \mu_n)=\mu$ is the $n$--tuple of eigenvalues
of $B$.

For each $s\in\mathbb N$, we separate $\mathcal H_n^s(\mathbb C)$
into two parts, i.e. $\mathcal H_n^s(\mathbb C)=\mathcal
H_n^{sr}\oplus \mathcal H_n^{sn}$, where $\mathcal H_n^{sr}$
(resp. $\mathcal H_n^{sn}$) consists of vector--valued resonant
(resp. nonresonant) homogeneous polynomials of degree $s$.
Correspondingly,  we separate the right hand side of \eqref{2.2}
into
$[f]_s+[\phi]_s-g_s=([f]_{sr}+[\phi]_{sr}-g_{sr})+([f]_{sn}+[\phi]_{sn}-g_{sn})\in\mathcal
H_n^{sr}\oplus \mathcal H_n^{sn}$. Since the operator $\mathcal
L_{B,s}$ is linear, equation \eqref{2.2} can be written in two
equations
\begin{eqnarray}
\mathcal L_{B,s}\phi_{sr}=[f]_{sr}+[\phi]_{sr}-g_{sr},\label{2.21}\\
\mathcal L_{B,s}\phi_{sn}=[f]_{sn}+[\phi]_{sn}-g_{sn},\label{2.22}
\end{eqnarray}
with $\phi_s=\phi_{sr}+\phi_{sn}\in\mathcal H_n^{sr}\oplus
\mathcal H_n^{sn}$. For equation \eqref{2.21}, we choose
$g_{sr}=[f]_{sr}+[\phi]_{sr}$, and consequently it has the trivial
solution $\phi_{sr}=0$. Since $\mathcal L_{B,s}$ is invertible on
$\mathcal H_n^{sn}$, equation \eqref{2.22} has always a unique
solution $\phi_{sn}$ for any $g_{sn}\in\mathcal H_n^{sn}$. We
choose the solution of \eqref{2.22} with $g_{sn}=0$. Then
$\phi_s(y)=\phi_{sn}(y)$ contains only nonresonant monomials,
and $g_s(y)=g_{sr}$ contains only resonant monomials.

The above proof shows that the normalization $\Phi(y)=y+\phi(y)$
has its nonlinear part consisting of nonresonant terms, i.e. all
monomials $y^me_i$ in the $i$th component satisfying
$\mu^m-\mu_i\ne 0$. The distinguished normal form $G(y)=By+g(y)$
has its nonlinear part consisting of resonant terms. Moreover, we
know from the above proof that the distinguished normal form and
normalization are both unique. This proves the lemma.

The next result characterizes the first integral of the
distinguished normal form for the given diffeomorphism $F(x)$.

\begin{lemma}\label{l22}
Let $G(y)=By+g(y)$ be the distinguished normal form of
$F(x)=Bx+f(x)$ via the distinguished normalization
$x=\Phi(y)=y+\phi(y)$. If $F(x)$ has an analytic first integral,
then $G(y)$ has a first integral either analytic or formal with
its nonlinear term all resonant.
\end{lemma}

We should mention the difference on resonance between first
integrals and normal forms including normalization. By definition
a monomial $y^m$ in a first integral of $G(y)$ is {\it resonant}
if $\mu^m=1$.


\noindent{\it Proof}. By the assumption $F$ and $G$ are conjugate, i.e. $F\circ \Phi=\Phi\circ
G$. If $V(x)$ is an analytic first integral of $F(x)$, then
$W(y)=V\circ \Phi(y)$ is an analytic or formal first integral of
$G(y)$ because $W(G(y))=W\circ G(y)=V\circ \Phi\circ G(y)=V\circ
F\circ \Phi(y)=V\circ \Phi(y)=W(y)$, where we have used the fact
that $V(F(x))=V(x)$ for all $x\in(\mathbb C^n,0)$.

Next we prove that the first integral $W$ of $G$ consists of
resonant monomials, i.e. its each monomial $y^m$ modulo
coefficient satisfying $\mu^m=1$. Indeed, $W$ is a first integral
of $G(y)$ means that $W(G(y))=W(y)$ for all $y\in(\mathbb C^n,0)$.
We rewrite this last equation as
\begin{equation}\label{2.3}
W(By)-W(y)=W(By)-W(By+g(y)).
\end{equation}
From the above proof we can set
$W(y)=\sum\limits_{s=r}\limits^\infty W_s(y)$ with $W_s(y)$
homogeneous polynomial of degree $s$ in $y$ and $W_r(y)\not\equiv
0$. Re--expanding $W(By+g(y))$ in Taylor series in $y$, we get
from \eqref{2.3} that
\begin{equation}\label{2.4}
W_s(By)-W_s(y)=R_s(y), \quad s=r,r+1,\ldots,
\end{equation}
where $R_s(y)$ is inductively known and $R_r(y)=0$. From Lemma
\ref{l20} we get that the linear operator $\mathcal L_s:
H^s(\mathbb C^n)\rightarrow H^s(\mathbb C^n)$ defined by $\mathcal
L_s\phi(y)=\phi(By)-\phi(y)$ for $\phi\in H^s(\mathbb C^n)$ has
the spectrum $\{\mu^m-1;\,\,m\in\mathbb Z_+^n,|m|\ge 2\}$, where $
H^s(\mathbb C^n)$ is the linear space of scalar complex
homogeneous polynomials of degree $s$ in $n$ variables. So
equation \eqref{2.4} with $s=r$, i.e. $\mathcal L_rW_r(y)=0$, has
only the resonant homogeneous polynomial solution. So $W_r$ must
be a resonant homogeneous polynomial of degree $r$. For equation
\eqref{2.4} with $s>r$, it is easy to know from the right hand
side of \eqref{2.3} that the $R_s(y)$ is constructed from
$W_r,\ldots,W_{s-1}$. By induction we assume that
$W_r,\ldots,W_{s-1}$ are all resonant. So in order for proving
$R_s(y)$ to be resonant, we only need to prove that each term
$z^m$ in $W(By+g(y))$ with $m\in\mathbb Z_+^n$, $|m|<s$ and
$z=By+g(y)$ contains only resonant monomials in $y$. For this aim,
we set $B$ in lower triangular Jordan normal form with
\[
B=\left(\begin{array}{ccccc} \mu_1 & 0 & 0 & 0 & 0\\
\sigma_1 & \mu_2 & 0 & 0 & 0 \\
0 & \sigma_2 & \ddots & 0 & 0 \\
0 & 0 & \sigma_{n-2} & \mu_{n-1} & 0 \\
0 & 0 & 0 & \sigma_{n-1} & \mu_n \end{array}\right),
\]
where $\sigma_i=0$ or $1$ for $i=1,\ldots,n-1$, and if
$\sigma_i=1$ then $\mu_{i-1}=\mu_i$. In these notations we have
$z^m=\prod\limits_{i=1}\limits^{n}(\sigma_{i-1}y_{i-1}+\mu_iy_i+g_i(y))^{m_i}$,
where $\sigma_0=0$. Using the binormal expansion we get
$(\sigma_{i-1}y_{i-1}+\mu_iy_i+g_i(y))^{m_i}=\sum\limits_{j=0}\limits^{m_i}\left(\begin{array}{c}
m_i\\
j\end{array}\right)(\sigma_{i-1}y_{i-1}+\mu_iy_i)^{m_i-j}(g_i(y))^{j}$.
Since $g_i(y)$ contains only resonant monomials, i.e. its each
monomial $y^k$ satisfies $\mu^k=\mu_i$, it follows that each
monomial $y^l$ in $(g_i(y))^{j}$ satisfies $\mu^l=\mu_i^j$. In
addition, if $\sigma_{i-1}\ne 0$ then $y_{i-1}$ and $y_i$ satisfy
the same resonant conditions because the eigenvalues corresponding
to $y_{i-1}$ and $y_i$ are the same. These two facts show that the
monomial $y^r$ in
$(\sigma_{i-1}y_{i-1}+\mu_iy_i)^{m_i-j}(g_i(y))^{j}$ verifies
$\mu^r=\mu_i^{m_i}$. Consequently each monomial $y^q$ in $z^m$
satisfies $\mu^q=\mu_1^{m_1}\ldots\mu_n^{m_n}=1$, where we have
used the fact that each monomial $z^m$ of $W_i(z)$ for $l\le i<s$
is resonant. This proves the lemma.

We now study the number of functionally independent analytic or
formal first integrals for analytic diffeomorphisms.

\begin{lemma}\label{l23}
The analytic diffeomorphism $F(x)=Bx+f(x)$ has at most $d_\mu$,
the rank of resonant set $\mathcal D$ of $B$, functionally
independent analytic or formal first integrals.
\end{lemma}

\noindent{\it Proof}. Let $G(y)=By+g(y)$ be the distinguished normal form of $F(x)$
through a normalization tangent to the identity, and let $W(y)$ be
an analytic or a formal first integral of $G$. It follows from
Lemma \ref{l22} that $W(y)$ is also a first integral of the linear
map $\mu y$ with $\mu y=(\mu_1 y_1,\ldots,\mu_n y_n)$, because $
W(y)$ consists of resonant monomials and each monomial $y^m$ in
$W$ satisfies $(\mu y)^m=\mu^my^m=y^m$.

Since $F(x)$ and $G(y)$ are conjugate via an analytic or a formal
normalization tangent to the identity, they have the same number
of functionally independent analytic or formal first integrals. So
we only need to prove that $G(y)$ has at most $d_{\mu}$
functionally independent analytic or formal first integrals. The
above proof shows that the the number of functionally independent
first integrals of $G(y)$ does not exceed that of $\mu y$. We now
turn to prove that the linear diffeomorphism $\mu y$ has exactly
$d_\mu$ functionally independent analytic  first integrals.

Let $m_1,\ldots,m_{d_\mu}$ be the linearly independent elements of
$\mathcal D$, and each vector $m_j$ for $j=1,\ldots,d_{\mu}$ is
simple, i.e. it cannot be divided by a positive integer no less
than $2$. Then $y^{m_1},\ldots,y^{m_{d_\mu}}$ are the $d_\mu$
functionally independent analytic first integrals of $\mu y$,
because $(\mu y)^{m_k}=y^{m_k}$ for $k=1,\ldots,d_\mu$ and their
gradients $\nabla
(y^{m_k})=((m_{k1}/y_1,\ldots,m_{kn}/y_n)y^{m_k}$,
$k=1,\ldots,d_\mu$, are linearly independent in an open dense
subset of $\mathbb C^n$. For any $m^*\in\mathcal D$, it linearly
depends on $m_1,\ldots,m_{d_\mu}$, and so $y^{m^*}$ functionally
depends on $y^{m_1},\ldots,y^{m_{d_\mu}}$. Also the proof of Lemma
\ref{l22} shows that any formal or analytic first integral of $\mu
y$ consists of resonant monomials.  This implies that any analytic
or formal first integral of $\mu y$ is an analytic or a formal
function of $y^{m}$ with $m\in\mathcal D$ simple. We have proved
that $\mu y$ has exactly $d_\mu$ functionally independent analytic
or formal first integrals. Consequently $F(x)$ has  at most
$d_\mu$ functionally independent analytic or formal first
integrals. The proof is completed.

We now study the expression of the normal form $G(y)$ and  the
relation between the orbits of $G(y)=By+g(y)$ and of $\mu y$.
\begin{lemma}\label{l24}
Assume that the analytic diffeomorphsim $F(x)=Bx+f(x)$ has $n-1$
functionally independent analytic first integrals, and that
$G(y)=By+g(y)$ is the distinguished normal form of $F(x)$. If $B$
has at least one eigenvalue  not equal to one in modulus, then the
following statements hold.
\begin{itemize}

\item[$(a)$] The resonant set of $B$ has the rank $d_\mu=n-1$.

\item[$(b)$] $B$ is diagonal, and
$G(y)=(\mu_1y_1(1+p_1(y)),\mu_2y_2(1+p_2(y),\ldots,\mu_ny_n(1+p_n(y))$
with $p_i(0)=0$ for $i=1,\ldots,n$. Moreover,
$p_1(y),\ldots,p_n(y)$ satisfy the functional equations
\[
(1+p_1(y))^{m_{k1}}\ldots (1+p_n(y))^{m_{kn}}=1,\quad
k=1,\ldots,n-1,
\]
where $(m_{k1},\ldots,m_{kn})=m_k\in\mathcal D$, $k=1,\ldots,n-1$,
are linearly independent and simple.

\item[$(c)$] The generic orbits of both $G(y)$ and $\mu y$ are contained in
the same orbits of some vector field.
\end{itemize}
\end{lemma}


\noindent{\it Proof}. $(a)$ Since the eigenvalues $\mu$ have modulus not all equal to 1,
it implies that the rank $d_\mu$ of the resonant set of $B$ is
less than or equal to $n-1$. By the assumption and Lemma \ref{l23}
we get that $d_\mu=n-1$.

\noindent $(b)$ Let $V_1(x),\ldots,V_{n-1}(x)$ be the $n-1$
functionally independent analytic first integrals without
constants of $F(x)$. Then it follows from the proof of Lemma
\ref{l22} that $W_i(y)=V_i\circ \Phi(y)$, $i=1,\ldots,n-1$, are
the functionally independent analytic or formal first integrals of
$G(y)$, where $x=\Phi(y)$ is the distinguished normalization from
$F(x)$ to $G(y)$. By Ziglin's lemma \cite{Zi} (see also the
appendix of \cite{It}), we can assume without loss of generality
that the lowest order parts $W_i^0(y)$ of $W_i(y)$ for
$i=1,\ldots,n-1$ are functionally independent. Otherwise it can be
done by polynomial combination of these integrals with complex
coefficients.

Since $W_i(y)$ consists of resonant monomials, we get that
$W_i(By+g(y))=W_i(y)$ and $W_i(\mu y)=W_i(y)$ for all $y$ in
$(\mathbb C^n,0)$. Equating the lowest order terms in $y$ of these
last two equations gives $W_i^0(By)=W_i^0(y)=W_i^0(\mu y)$. Set
$B=U+N$ with $U=\mbox{diag}(\mu_1,\ldots,\mu_n)$ and $N$ in the
nilpotent lower triangular normal form, i.e. we have
$Ny=(0,\sigma_1y_1,\ldots,\sigma_{n-1}y_{n-1})$. Then we get from
$W_i^0(\mu y+Ny)=W_i^0(\mu y)$ that
\begin{equation}\label{2.10}
\langle\nabla W_i^0(\mu y+\theta_y Ny),\,Ny\rangle=0,\quad
i=1,\ldots,n-1,
\end{equation}
where $\theta_y\in (0,1)$ and $\mu y=Uy$. Since
$W_1^0,\ldots,W_{n-1}^0$ are functionally independent in $(\mathbb
C^n,0)$, we can assume without loss of generality that
\[
\Delta^*(y)=\det\left(\begin{array}{ccc}
\frac{\partial W_1^0}{\partial x_2}(z_y) & \cdots & \frac{\partial W_1^0}{\partial x_n}(z_y)\\
\vdots & \ddots & \vdots\\
\frac{\partial W_{n-1}^0}{\partial x_2}(z_y) & \cdots &
\frac{\partial W_{n-1}^0}{\partial x_n}(z_y)
\end{array}\right)\ne 0,
\]
in an open subset of $(\mathbb C^n,0)$, where $z_y=\mu y+\theta_y
Ny$. Otherwise it can be done by rearranging the order of the
coordinates, and meanwhile the Jordan normal form $B$ keeps in the
same form. Hence equation \eqref{2.10} has the unique solution
$Ny=0$, and consequently $N=0$. This proves that $B$ is diagonal.

For proving $G(y)$ to have the special type of normal form,
instead of the $n-1$ functionally independent first integrals
$W_i(y)$ we consider the $n-1$ functionally independent monomial first
integrals $H_k(y)=y^{m_k}$ for $k=1,\ldots,n-1$ with
$m_1,\ldots,m_{n-1}\in\mathcal D$ being linearly independent and
simple. Here the existence of the $n-1$ monomial first integrals $H_k(y)$ follows from the facts that since
$W_1(y),\ldots,W_{n-1}(y)$ are functional independent, and so the Inverse Function Theorem implies that there exist
$n-1$ functionally independent monomials $H_1(y),\ldots,H_{n-1}(y)$ in one--to--one way on an open and dense subset such that they are first integrals of $G(y)$, i.e. $H_i(G(y))=H_i(y)$ hold in an open and dense subset of $(\mathbb C^n,0)$ for $i=1,\ldots,n-1$. Now $H_i(y)$ are monomials and $G(y)$ is an analytic function or a formal series force that  $H_i(G(y))=H_i(y)$ must hold in $(\mathbb C^n,0)$, because by expanding these last equations and equating the homogeneous terms of the same order, we get a series of homogeneous polynomial equations. They hold in an open and dense subset of $(\mathbb C^n,0)$ and so must hold in $(\mathbb C^n,0)$.

Since $B$ is diagonal, we have
$G(y)=(\mu_1y_1+g_1(y),\ldots,\mu_ny_n+g_n(y))$. In the next proof
we distinguish two cases: either for any $k\in\{1,\ldots,n-1\}$,
$H_k(y)$ does not contain $y_1$; or there exists some
$k_0\in\{1,\ldots,n-1\}$ for which $H_{k_0}(y)$ contains $y_1$.

In the former, the first components of $m_k$ for $k=1,\ldots,n-1$
are all zero. Hence we have $\langle \overline m_k,\overline
\mu\rangle=0$ for $k=1,\ldots,n-1$, where $\overline
m_k=(m_{k2},\ldots,m_{kn})$ and $\overline
\mu=(\mu_2,\ldots,\mu_n)$. Obviously, $\overline
m_1,\ldots,\overline m_{n-1}$ are linearly independent. In $n-1$
dimensional case, $\overline \mu$ satisfies $n-1$ linearly
independent resonant relations, it follows from the proof of Lemma
\ref{l23} that $\overline \mu$ has its components all having
modulus $1$. By the assumption we must have $|\mu_1|\ne 1$. Since
$G(y)$ is in normal form, any nonlinear monomial $y^k$ in the
first component of $G(y)$ satisfies $\mu_1=\mu^k$. So we have
$|\mu_1|=|\mu_1|^{k_1}$, i.e. $k_1=1$. This proves that $y_1$
divides $g_1(y)$.

In the latter, set $g_1(y)=y_1p_1(y)+q_1(y)$, where
$p_1(y)=O(|y|)$, and $q_1(y)=O(|y|^2)$ is independent of $y_1$.
Using the fact that $H_{k_0}(y)$ are the first integrals of both
$G(y)$ and $\mu y$, i.e. $H_{k_0}(G(y))=H_{k_0}(y)=H_{k_0}(\mu
y)$, we obtain that
\[
\begin{array}{l}
(\mu_1y_1+y_1p_1(y)+q_1(y))^{m_{k_01}}(\mu_2y_2+g_2(y))^{m_{k_02}}
\ldots (\mu_2y_2+g_2(y))^{m_{k_0n}} \\
\quad
=(\mu_1y_1)^{m_{k_01}}(\mu_2y_2)^{m_{k_02}}\ldots(\mu_ny_n)^{m_{k_0n}}.
\end{array}
\]
In this last equation by setting $y_1=0$  gives
\[
(q_1(y))^{m_{k_01}}(\mu_2y_2+g_2(0,y_2,\ldots,y_n))^{m_{k_02}}
\ldots (\mu_2y_2+g_2(0,y_2,\ldots,y_n))^{m_{k_0n}}\equiv 0.
\]
This verifies that $q_1(y)\equiv 0$, because $m_{k_01}\ne 0$ and
$\mu_i\ne 0$ for $i=2,\ldots,n$.

The proof of the above two cases shows that the first component of
the normal form $G(y)$ is of the form $\mu_1y_1(1+p_1(y))$.
Working out in the same line we can prove that the $j$th component
of $G(y)$ is of the form $\mu_jy_j(1+p_j(y))$ for $j=2,\ldots,n$.

Finally using the first integrals $H_k(y)$ of $G(y)$ and of $\mu
y$, we obtain that
\[
(\mu_1y_1(1+p_1(y)))^{m_{k1}}\ldots
(\mu_ny_n(1+p_n(y)))^{m_{kn}}=(\mu y)^{m_k},\quad k=1,\ldots,n-1.
\]
Simplifying these last equations yields
\[
(1+p_1(y))^{m_{k1}}\ldots (1+p_n(y))^{m_{kn}}=1,\quad
k=1,\ldots,n-1,
\]
This proves statement $(b)$.

\noindent $(c)$ We will use the notations given in the proof of
statement $(b)$. The above proof shows that
$W_1(y),\ldots,W_{n-1}(y)$ are functionally independent first
integrals of both $G(y)$ and $\mu y$. So the level surfaces of
$W_i$ for $i=1,\ldots,n-1$ are invariant under the action of
either $G(y)$ or $\mu y$. This implies that each orbit of
$G(y)$ and of $\mu y$ is contained in the level surfaces of $W_i$
for $i=1,\ldots,n-1$ and so in their intersection. Clearly the
intersection is one dimensional in the full Lebesgue measure subset of $(\mathbb C^n,0)$ because of the functionally
independence of the $n-1$ first integrals.

Define a vector field in $(\mathbb C^n,0)$ by
\[
\mathcal Z(y)=\nabla W_1(y)\times\ldots\times\nabla
W_{n-1}(y),\quad \mbox{ for } y\in(\mathbb C^n,0),
\]
where $\times$ denotes the cross product of vectors in $\mathbb
C^n$. Recall that the cross product of $n-1$ vectors in $\mathbb
C^n$, saying $v_1,\ldots,v_{n-1}$, is again a vector, and is
defined by
\[
\langle v_1\times\ldots\times v_{n-1},
w\rangle=\det\left(\begin{array}{c}w\\v_1\\ \vdots\\
v_{n-1}\end{array}\right),
\]
for arbitrary $w\in \mathbb C^n$. By the very definition of the
cross product, it is easy to check that $W_k$ for $k=1,\ldots,n-1$
are first integrals of the vector field $\mathcal Z(y)$. So the
orbits of the vector field $\mathcal Z(y)$ are contained in the
intersections of the level surfaces of $W_1(y),\ldots,W_{n-1}(y)$.
This proves that both orbits of the diffeomorphisms $G(y)$ and
$\mu y$ starting at the same generic point are contained in the same
orbit of $\mathcal Z(y)$. Recall that the generic points are those ones which are located in a full Lebesgue measure subset of $(\mathbb C^n,0)$. We finish the proof of the lemma.

Next we prove that the nonresonant spectrum of a linear operator
related to $B$ is bounded from below in modulus.
\begin{lemma}\label{l25}
Assume that the diffeomorphism $F(x)=Bx+f(x)$ has $n-1$
functionally independent analytic first integrals and that $B$ has
at least one eigenvalue with modulus not equal to $1$. Then there
exits a $\sigma>0$ such that if $\mu^m-\mu_i\ne 0$ for
$m\in\mathbb Z_+^n$, $|m|\ge 2$ and $i=1,\ldots,n$, we have
$\left|\mu^m-\mu_i\right|\ge \sigma$.
\end{lemma}

\noindent{\it Proof}. By the assumption of the lemma we get from Lemma \ref{l24} that
there exist $n-1$ linearly independent vectors
$k_i=(k_{i1},\ldots,k_{in})\in \mathbb Z_+^n$ such that
$\mu^{k_i}=1$ for $i=1,\ldots,n-1$. This follows that
\begin{equation}\label{2.5}
k_{i1}\log|\mu_1|+\ldots +k_{in}\log|\mu_n|=0,\quad
i=1,\ldots,n-1.
\end{equation}
Since $k_1,\ldots,k_{n-1}$ are linearly independent, we can assume
without loss of generality that
\[
\Delta=\det\left(\begin{array}{ccc}k_{11} & \cdots & k_{1,n-1}\\
\vdots & \ddots & \vdots \\
k_{n-1,1} & \cdots & k_{n-1,n-1}\end{array} \right)\ne 0.
\]
Then we get from equation \eqref{2.5} using the Cram's rule that
\begin{equation}\label{2.6}
\log|\mu_j|=\frac{\delta_j}{\Delta}\,\log|\mu_n|,\quad j=1,\ldots,
n-1,
\end{equation}
where $\delta_j$'s$,\Delta\in\mathbb Z$. Moreover we have $\log
|\mu_n|\ne 0$. Otherwise all $\mu_j$ have modulus 1, a
contradiction with the assumption of the lemma. Using \eqref{2.6}
we get that for any $m\in\mathbb Z_+^n$ and $j=1,\ldots,n$
\begin{eqnarray*}
\left|\mu^m-\mu_j\right|&\ge
&\left||\mu_1|^{m_1}\ldots|\mu_n|^{m_n}-|\mu_j|\right|\\
 &=&
 |\mu_n|^{\frac{\delta_j}{\Delta}}\left||\mu_n|^{\frac{m_1\delta_1+\ldots+(m_j-1)\delta_j+\ldots+m_n\delta_n}{\Delta}}-1\right|,
\end{eqnarray*}
where $\delta_n=\Delta$. Set $\alpha=|\mu_n|^{1/\Delta}$. Since
$s_m=m_1\delta_1+\ldots+(m_j-1)\delta_j+\ldots+m_n\delta_n\in\mathbb
Z$, we have either $\alpha^{s_m}-1=0$ if $s_m=0$, or
$|\alpha^{s_m}-1|\ge \min\{|\alpha-1|, |\alpha^{-1}-1|\}\ne 0$ if
$s_m\ne 0$.

For $s_m\ne 0$, set
\[
\sigma_1=\min\{\alpha^{\delta_j}|\alpha-1|,
\alpha^{\delta_j}|\alpha^{-1}-1|;\, j=1,\ldots,n\},
\]
we have $\left|\mu^m-\mu_j\right|\ge \sigma_1$ for $j=1,\ldots,n$
and all $m\in\mathbb Z_+^n$ with $|m|\ge 2$ such that $s_m\ne 0$.

We now consider those $m$ such that $s_m=0$. From $\mu^{k_j}=1$
for $j=1,\ldots,n-1$ we get that
\begin{equation}\label{2.7}
 k_{j1}\log \mu_1+\ldots+k_{jn}\log \mu_n=\log 1,\quad
j=1,\ldots,n-1,
\end{equation}
where $\log 1=2n\pi\sqrt{-1}$, $n\in\mathbb Z$, and the logarithms
are taken for complex numbers because the eigenvalues $\mu$ may be
complex. Solving \eqref{2.7} by the Cram's rule gives
\begin{equation}\label{2.8}
\log
\mu_j=\frac{2n_j\rho_j\pi\sqrt{-1}+\delta_j\log\mu_n}{\Delta},\quad
j=1,\ldots,n-1,
\end{equation}
where $\delta_j$ is the same as that of \eqref{2.6},
$\rho_j\in\mathbb Z$ is uniquely determined by $k_i$ for
$i=1,\ldots,n-1$, and $n_j\in\mathbb Z$ come from the expression
of $\log 1$. For the $m\in\mathbb Z_+^n$ such that $s_m=0$, we
have
\begin{eqnarray*}
\left|\mu^m-\mu_j\right|&=&\left|e^{m_1\log \mu_1}\ldots
e^{m_n\log \mu_n}-e^{\log \mu_j}\right|\\
&=&\left|e^{\frac{2\sum\limits_{k=1}\limits^{n}m_kn_k\rho_k\pi\sqrt{-1}+\sum\limits_{k=1}\limits^{n}m_k\delta_k\log\mu_n}{\Delta}}-e^{\frac{2n_j\rho_j
\pi\sqrt{-1}+\delta_j\log \mu_n}{\Delta}}\right|
\\
&=&\left|e^{\frac{2\sum\limits_{k=1}\limits^{n}m_kn_k\rho_k\pi\sqrt{-1}}{\Delta}}\,\mu_n^{\frac{\sum\limits_{k=1}\limits^{n}m_k\delta_k}{\Delta}}
-e^{\frac{2n_j\rho_j \pi\sqrt{-1}}{\Delta}}\,
\mu_n^{\frac{\delta_j}{\Delta}}\right|
\\
&=&|\mu_n|^{\frac{\delta_j}{\Delta}}\left|e^{\frac{2\sum\limits_{k=1}\limits^{n}m_kn_k\rho_k\pi\sqrt{-1}}{\Delta}}-e^{\frac{2n_j\rho_j
\pi\sqrt{-1}}{\Delta}}\, \right|,
\end{eqnarray*}
where we have used the fact $s_m=0$ in the fourth equality, i.e.
$\sum\limits_{k=1}\limits^{n}m_k\delta_k=\delta_j$. Since $\Delta,
\rho_k,\rho_j$ are given integers which are uniquely determined by
the linearly independent vectors $k_i$ for $i=1,\ldots,n-1$, by
the periodic property of the exponential functions with respect to
their pure imaginary parts, it follows that
\begin{equation}\label{2.9}
e^{\frac{2\sum\limits_{k=1}\limits^{n}m_kn_k\rho_k\pi\sqrt{-1}}{\Delta}}\quad
\mbox{ and } \quad e^{\frac{2n_j\rho_j \pi\sqrt{-1}}{\Delta}}
\end{equation}
both take only finitely many values for all possible choice of
$m_k,n_k,n_j$. Taking $\gamma$ to be the minimum of the modulus of
all possible differences of the two elements given in \eqref{2.9}.
By the assumption $|\mu^m-\mu_j|\ne 0$ we have $\gamma\ne 0$, that
is, the modulus of their difference has a nonzero minimum. Set
\[
\sigma_2=\min\{\alpha^{\delta_j}\gamma;\, j=1,\ldots,n\}.
\]
Then $\sigma=\min\{\sigma_1,\sigma_2\}$ is the data satisfying the
lemma, i.e. we have $|\mu^m-\mu_j|\ge \sigma$ for all $m\in\mathbb
Z_+^n$ with $|m|\ge 2$, $j=1,\ldots,n$ and $|\mu^m-\mu_j|\ne 0$.
We complete the proof of the lemma.

The last step is to prove that the normalization from $F(x)$ to
$G(y)$ is convergent. Lemma \ref{l25} shows that in the
analytic integrable case there does not appear small divisor
conditions. A folklore says that if no small divisor conditions
appear, it is convergent that the distinguished normalization
tangent to identity from a given analytic vector field or analytic
diffeomorphism to its distinguished normal form. In fact, it is
not the case. See the following example, the planar analytic
vector field $\dot x=x+\varphi(x,y)$, $\dot y=-y+\psi(x,y)$ is
always formally equivalent to $\dot x=xf(xy)$, $\dot y=-yg(xy)$.
The eigenvalues $\lambda_1=1,\lambda_2=-1$ of its linear part do
not satisfy small divisor conditions, because $0\ne
|q_1\lambda_1+q_2\lambda_2-1|\ge 1$ for $q_1,q_2\in \mathbb Z_+$.
But generally no results guarantee the convergence of the
normalization except for $f(xy)=g(xy)$ (see example 2.3 and the
remark following Theorem 2.4 of \cite{St09}).

The following result shows that for analytic integrable
diffeomorphism, the distinguished normalization is analytic.

\begin{lemma}\label{l26}
If the analytic diffeomorphism $F(x)=Bx+f(x)$ is analytic
integrable and $B$ has at least one eigenvalue not equal to $1$ in
modulus, then it is analytically conjugate to its normal form of
the type
\[
G(y)=(\mu_1y_1(1+p_1(y)),\ldots,\mu_ny_n(1+p_n(y))),
\]
where $p_1(y),\ldots,p_n(y)$ can be represented in analytic
functions of a single analytic function, and
$p_1(0)=\ldots=p_n(0)=0$.
\end{lemma}

\noindent{\it Proof}. Lemma \ref{l24} has showed that $F(x)$ has the prescribed normal
form, and that there exist $n-1$ linearly independent simple
resonant lattice $m_k\in\mathcal D$, $k=1,\ldots,n-1$, such that
\begin{equation}\label{l13} (1+p_1(y))^{m_{k1}}\ldots
(1+p_n(y))^{m_{kn}}=1,\quad k=1,\ldots,n-1.
\end{equation}
Since $m_1,\ldots,m_{n-1}$ are linearly independent, solving
\eqref{l13} yields that there exists an $\iota\in\{1,\ldots,n\}$
such that
$1+p_1(y),\ldots,1+p_{\iota-1}(y),1+p_{\iota+1}(y),\ldots,1+p_n(y)$
can be represented in functions of $1+p_\iota(y)$. More precisely,
for $j\in\{1,\ldots,\iota-1,\iota+1,\ldots,n\}$ we have
$1+p_j(y)=(1+p_\iota(y))^{p_j/q}$ with $p_j,q\in\mathbb Z$
uniquely determined by $m_k$'s. Obviously, $p_j$ is an analytic
function in $p_\iota$ if $|p_\iota|<1$. Of course, if $p_\iota(y)$
is locally analytic in $(\mathbb C^n,0)$, then
$p_1(y),\ldots,p_{\iota-1}(y),p_{\iota+1}(y),\ldots,p_n(y)$ will
be locally analytic in $(\mathbb C^n,0)$.

In what follows we assume without loss of generality that $l=1$.
From the proof of Lemma \ref{l21} the diffeomorphism
$F(x)=Bx+f(x)$ is formally transformed to $G(y)=By+g(y)$ by a
formal distinguished normalization $x=y+\phi(y)$. Set
\[
f_s(x)=\sum\limits_{m\in\mathbb Z_+^n,|m|\ge l}f_s^mx^m,\,\,\,\,
g_s(y)=\sum\limits_{m\in\mathbb Z_+^n,|m|\ge l}g_s^my^m,\,\,\,\,
\phi_s(y)=\sum\limits_{m\in\mathbb Z_+^n,|m|\ge l}\phi_s^my^m,
\]
for $s=1,\ldots,n$, where $f_s^m$, $g_s^m$ and $\phi_s^m$ are the
coefficients of $x^m$ and $y^m$ respectively, $l$ is the degree of
the lowest order term of $f(y)$, and $h_s$ is the $s$th component
of $h\in\{f,g,\phi\}$. Then by Lemmas \ref{l21} and \ref{l24}, and
comparing the coefficients of $y^m$ in the $k$th component of
\eqref{2.1}, we get that
\begin{equation}\label{2.14}
(\mu^m-\mu_k)\phi_k^m=\left[f_k\right]^m-\sum\limits_{r\in\mathbb
Z_+^m,\,r\precneqq m}\phi_k^r\mu^r P_r^{m-r}-\mu_k p_k^{m-e_k},
\end{equation}
where $[f_k]^m$ is the coefficient of $y^m$ in the expansion of
$f_k(y+\phi(y))$, $r\precneqq m$ means that $r\ne m$ and
$m_s-r_s\ge 0$ for $s=1,\ldots,n$, and $P_r^{m-r}$ is the
coefficient of $y^{m-r}$ of $P_r(y)=(1+p_1(y))^{r_1}\ldots
(1+p_n(y))^{r_n}$. Here we have used the fact that
$\phi_k(By+g(y))=\phi_k(\mu_1y_1(1+p_1(y)),\ldots,\mu_ny_n(1+p_n(y)))=\sum\limits_{m\in
\mathbb Z_+^n,\,|m|\ge l}\phi_k^m\mu^my^m(1+p_1(y))^{m_1}\ldots
(1+p_n(y))^{m_n}$.

For $m\in\mathbb Z_+^n$ such that $\mu^m=\mu_k$, we have from
\eqref{2.14} that
\begin{equation}\label{2.15}
\phi_k^m=0,\quad p_k^{m-e_k}=\mu_k^{-1}\left[f_k\right]^m,
\end{equation}
where we have used the fact that $\sum\limits_{r\in\mathbb
Z_+^m,\,r\precneqq m}\phi_k^r\mu^r P_r^{m-r}=0$. Because $P_r(y)$
contains only resonant term, it follows that $\mu^{m-r}=1$ if
$P_r^{m-r}\ne 0$, and so $\mu^r=\mu^m=\mu_k$. This implies that
the monomial $\phi_k^ry^r$ in $\phi_k(y)$ is resonant and so
$\phi_k^r=0$. Recall that since $F(x)=Bx+f(x)$ is a
diffeomorphism,  we have $\mu_k\ne 0$ for $k=1,\ldots,n$. Then we
have the estimation
\begin{equation}\label{2.151}
|p_k^{m-e_k}|\le \nu|\left[f_k\right]^m|,
\end{equation}
where $\nu=\max\{1/|\mu_k|;\,k=1,\ldots,n\}$.

For $m\in\mathbb Z_+^n$ such that $\mu^m\ne \mu_k$, we have from
\eqref{2.14} that
\begin{equation}\label{2.16}
p_k^{m-e_k}=0,\quad
\phi_k^m=\frac{1}{\mu^m-\mu_k}\left(\left[f_k\right]^m-\sum\limits_{r\in\mathbb
Z_+^m,\,r\precneqq m}\phi_k^r\mu^r P_r^{m-r}\right).
\end{equation}
Furthermore, by Lemma \ref{l25} and the fact that $\mu^{m-r}=1$ if
$P_r^{m-r}\ne 0$ we have the estimation
\begin{eqnarray*}
\left|\frac{1}{\mu^m-\mu_k}\sum\limits_{r\in\mathbb
Z_+^m,\,r\precneqq m}\phi_k^r\mu^r
P_r^{m-r}\right|&=&\sum\limits_{r\in\mathbb Z_+^m,\,r\precneqq
m}\left(1+\frac{|\mu_k|}{|\mu^m-\mu_k|}\right)|\phi_k^r||P_r^{m-r}|\nonumber\\
&\le &\delta\sum\limits_{r\in\mathbb Z_+^m,\,r\precneqq
m}|\phi_k^r||P_r^{m-r}|,
\end{eqnarray*}
where $\delta=1+\sigma^{-1}\max\{|\mu_k|;\, k=1,\ldots,n\}$.
Recall that $\sigma$ is the data given in Lemma \ref{l25}. Then we
have the estimation for $\phi_k^m$ given in \eqref{2.16}
\begin{equation}\label{2.161}
|\phi_k^m|\le \sigma^{-1}|[f_k]^m|+\delta\sum\limits_{r\in\mathbb
Z_+^n,\,r\precneqq m}|\phi_k^r||P_r^{m-r}|.
\end{equation}

Having the above estimations we can use the majorant series to
prove the convergence of $\phi_s(y)$ and of $p_s(y)$ for
$s=1,\ldots,n$. For a series $h_s(y)=\sum\limits_{m\in\mathbb
Z_+^n}h_s^my^m$, we define $\hat h_s(y)=\sum\limits_{m\in\mathbb
Z_+^n}|h_s^m|y^m$. For two scalar series $\xi(y)$ and $\eta(y)$,
we say that the latter is a majorant series of the former, denoted
by $\xi(y)\preccurlyeq \eta(y)$, if $\left|\xi^m\right|\le \eta^m$ and
$\eta^m\ge 0$, where $\xi^m$ and $\eta^m$ are the coefficients of
$y^m$ in the series $\xi(y)$ and $\eta(y)$, respectively. Under
this notation we have $h_s\preccurlyeq\hat h_s$ for the scalar
series $h_s$. We refer the readers to \cite{Hi} for more detail
information on the majorant series.

Since $F(x)=Bx+f(x)$ is analytic in $(\mathbb C^n,0)$, by the
Cauchy inequality there exists a polydisc
$\Omega_\rho=\{|x_s|<\rho;\, s=1,\ldots,n\}$ in which we have
\[
|f_s^m|\le M\rho^{-|m|},\qquad \mbox{ for } s=1,\ldots,n,
\]
where $M=\max\limits_s\sup\limits_{\partial \Omega_\rho}\{|f_s|\}$
and $f_s$ is the $s$th component of $f$. Clearly,
\[
\widetilde f(x)=\sum\limits_{m\in\mathbb Z_+^n}M\rho^{-|m|}x^m,
\]
is convergent in $\Omega_\rho$, and $\hat f_s(x)\preccurlyeq
\widetilde f(x)$ for $s=1,\ldots,n$. So the majorant series $\hat
f_s(x)$ of $f_s(x)$ is convergent in $\Omega_\rho$, and
consequently is analytic in the domain.

Since $\phi_s(y)$ and $p_s(y)$ have the coefficients satisfying
\eqref{2.15} and \eqref{2.16} with the estimates \eqref{2.151} and
\eqref{2.161}, by some calculations we get that
\begin{eqnarray}\label{2.19}
\sum\limits_{k=1}\limits^n(\phi_k(y)+p_k(y))&\preccurlyeq &
\sum\limits_{k=1}\limits^n(\hat\phi_k(y)+\hat
p_k(y))\\
&\preccurlyeq & n(\sigma^{-1}+\nu)\widetilde f(y+\hat
\phi(y))+\delta\sum\limits_{k=1}\limits^n\left(\hat\phi_k(y(1+\hat
p(y))-\hat\phi_k(y)\right),\nonumber
\end{eqnarray}
where $y(1+\hat p(y))=(y_1(1+\hat p_1(y)),\ldots,y_n(1+\hat
p_n(y)))$. For simplicity to notation we set
$\gamma=n(\sigma^{-1}+\nu)$. By the very definition of $\hat
\phi_s$ and $\hat p_s$, in order for proving the convergence of
$\sum\limits_{k=1}\limits^n(\hat\phi_k(y)+\hat p_k(y))$ in
$\Omega_{\rho_*}$ with $\rho_*\in (0,\rho)$ to be specified later
on, we only need to prove it when $y_1=\ldots=y_n=z$ and $|u|\le
\rho_*$. For doing so, we set
\[
\left. U(z)=\sum\limits_{k=1}\limits^n(\hat\phi_k(y)+\hat
p_k(y))\right|_{y_1=\ldots=y_n=z}.
\]
In fact, in $U(z)$ we can use only $\hat p_1(y)$, instead of
$\sum\limits_{k=1}\limits^n\hat p_k(y)$, because at the beginning
of the proof of this lemma we have proved that
$p_2(y),\ldots,p_n(y)$ can be represented in functions of
$p_1(y)$. Since the lowest order terms of $\hat \phi(y)$ and of
$\hat p_k(y)$ have degree no less than $1$, it follows that $U(z)$
must be divided by $z$. Set $U(z)=V(z)z$. We get from \eqref{2.19}
that
\begin{equation}\label{2.20}
V(z)\preccurlyeq \gamma z\widetilde f_*(V(z))+\delta \left((1+z
V(z))V(z(1+z V(z)))-V(z)\right),
\end{equation}
where $\widetilde f_*(V(z))$ is $\widetilde
f(z(1+V(z)),\ldots,z(1+V(z)))$ divided by $z^2$, and it is
analytic as a function of $V$. For obtaining \eqref{2.20} we have
used the facts that $\widetilde f(y+\hat \phi(y))\preccurlyeq
\widetilde f(z+W(z),\ldots,z+W(z))$ and that
\[
\sum\limits_{k=1}\limits^n\left(\hat\phi_k(y(1+\hat
p(y))-\hat\phi_k(y)\right)\preccurlyeq W(z(1+W(z)))-W(z).
\]

Set
\[
T(h,z)=h-\gamma z\widetilde f_*(h)-\delta \left((1+z h)h(z(1+z
h))-h\right).
\]
For studying the existence of analytic solution, saying $h(z)$, of
$T(h,z)=0$, we introduce an auxiliary function
\[
\Lambda(h,z)=h-\gamma z\widetilde f_*(h)-\delta \left((1+z
h)h-h\right).
\]
Obviously $\Lambda$ is analytic in $h$ and $z$, because
$\widetilde f_*$ is an analytic function in $h$. Some easy
calculations show that
\[
\Lambda(0,0)=0,\qquad \left.\frac{\partial \Lambda}{\partial
h}\right|_{(h,z)=(0,0)}=1.
\]
By the Implicit Function Theorem the equation $\Lambda(h,z)=0$ has
a unique analytic solution, denoted by $h_0(z)$, in a neighborhood
of $0$ in $C$.

Choose $\rho_1>0$ satisfying $\rho_1<\min\{1,\rho\}$ for which
$h_0(z)$ is analytic in $B_{\rho_1}(0)=\{z\in\mathbb C;\,
|z|<\rho_1\}$ and $\|h_0\|=\sup\{|h_0(z)|;\, z\in
B_{\rho_1}(0)\}<1$. Then the functional equation $T(h,z)=0$ has an
analytic solution $h(z)$ defined in $B_{\rho_1/3}(0)$. Comparing
\eqref{2.20} with $T(h,z)$, it follows that $h(z)$ is a majorant
series of $V(z)$. Hence $V(z)$ is analytic in $B_{\rho_1/3}(0)$,
and consequently $\sum\limits_{k=1}\limits^n(\hat \phi_k+\hat
p_k)$ is analytic in the ball. This proves that the distinguished
normalization from $F(x)=Bx+f(x)$ to $G(y)=By+g(y)$ is analytic,
that is, $F(x)$ is analytically conjugate to its distinguished
normal form. We complete the proof of the lemma.

Having the above preparations we can prove Theorem \ref{t2}.

\noindent{\bf Proof of Theorem \ref{t2}}: {\it Sufficiency}. By
the assumption of the theorem the monomials $H_k(y)=y^{m_k}$ for
$k=1,\ldots,n-1$ are $n-1$ functionally independent analytic first
integrals of $G(y)$, where $m_1,\ldots,m_{n-1}$ are the linearly
independent resonant lattices given in Theorem \ref{t2}. Let
$x=\Phi(y)$ be the analytic conjugation tangent to the identity
from $F(x)$ to $G(y)$, and let $y=\Psi(x)$ be its inverse. Then
$\Psi(x)$ is analytic, tangent to the identity and satisfies
$\Psi\circ F=G\circ \Psi$. On the other hand, using the conjugate
condition and $H_k\circ G(y)=H_k(y)$ for $y\in(\mathbb C^n,0)$ we
get that $H_k \circ \Psi\circ F(x)=H_k\circ G\circ
\Psi(x)=H_k\circ \Psi(x)$. This proves  that $H_k\circ \Psi(x)$,
$k=1,\ldots,n$, are $n-1$ analytic first integrals of $F(x)$.
Furthermore, by the functional independence of $H_1,\ldots,
H_{n-1}$ and $y=\Psi(x)$ tangent to identity, it follows easily
that $H_1\circ \Psi(x),\ldots,H_{n-1}\circ \Psi(x)$ are
functionally independent in $(\mathbb C^n,0)$. This proves that
$F(x)$ has $n-1$ functionally independent analytic first
integrals, and consequently is analytic integrable in
$(\mathbb C^n,0)$.

\noindent{\it Necessity}. The proof follows from Lemmas \ref{l24}
and \ref{l26}. We have completed the proof of the theorem. \hskip
5.6in $\Box$

\section{ Proof of Theorem \ref{t1}}

\setcounter{section}{3}
\setcounter{equation}{0}\setcounter{theorem}{0}

We should mention that the main idea of the proof follows from that of Theorem \ref{t2}.
Here we present a sketch proof and mainly concern the parts of the proof which are
different from those given in the proof of Theorem \ref{t2}.

\noindent {\it Sufficiency}. By the assumption $R_\lambda=n-1$,
there exist $n-1$ linearly independent vectors
$m_i=(m_{i1},\ldots,m_{in})\in\mathbb Z_+^n$, $i=1,\ldots,n$, such
that $\langle m_i,\,\lambda\rangle=0$. This implies that
$y^{m_i}$, $i=1,\ldots,n-1$ are $n-1$ functionally independent
analytic first integrals of \eqref{e2}. Let $y=\psi(x)$ be the
analytic transformation tangent to the origin from \eqref{e2} to
\eqref{e1} in a neighborhood of the origin. Then $\psi^{m_i}$,
$i=1,\ldots,n-1$, are the $n-1$ functionally independent analytic
first integrals of \eqref{e1}.

\noindent{\it Necessity}.  Denote by $\mathcal X$ the vector
fields induced by system \eqref{e1}. Set ${\mathcal X}={\mathcal
X}_1+{\mathcal X}_h$ with ${\mathcal X}_1$ and ${\mathcal X}_h$
the linear and higher order terms, respectively. Furthermore we
separate ${\mathcal X}_1={\mathcal X}_1^s+{\mathcal X}_1^n$ with
${\mathcal X}_1^s=\langle A_1x,\partial_x\rangle$ the {\it
semisimple part} and ${\mathcal X}_1^n=\langle
A_2x,\partial_x\rangle$ the {\it nilpotent part} of ${\mathcal
X}_1$ respectively, where $A=A_1+A_2$. Without loss of generality,
we can assume that
\[
{\mathcal
X}_1^s:=\sum\limits_{i=1}\limits^{n}\lambda_ix_i\frac{\partial}{\partial
x_i}.
\]
Recall that the vector field $\mathcal X$ is in {\it normal form}
is equivalent to that the Lie bracket of  ${\mathcal X}_1^s$ and
${\mathcal X}_h$ vanishes, i.e. $[{\mathcal X}_1^s, {\mathcal
X}_h]=0$.

For a given analytic differential system or vector field, by the
Poincar\'e-Dulac normal form theorem it can always be transformed
to its distinguished normal form by a distinguished normalization. Let
\begin{equation}\label{3.1}
\dot y= Ay+g(y),
\end{equation}
be the distinguished normal form of \eqref{e1} obtained from the
normalization $x=\Phi(y)=y+\varphi(y)$. Then the vector field
associated with \eqref{3.1} is $\mathcal
Y(y)=(D\Phi(y))^{-1}(A+f)\circ \Phi(y)$, where $D\Phi(y)$ denotes
the Jacobian matrix of $\Phi(y)$.

We claim that if $H(x)$ is an analytic first integral of
\eqref{e1}, then $V(y)=H(y+\varphi(y))$ is an analytic or formal
first integral of \eqref{3.1}, and all its monomials are resonant.
The proof is similar to Lemma 2.3 of \cite{XZ08}. The difference
is that now $A$ is priori not necessary diagonal. We now prove the
claim. That $H(x)$ is an analytic first integral of \eqref{e1} is
equivalent to $ \langle \partial_x H,\,Ax+F(x)\rangle=0$.  By the
chain rule, it follows that $\langle \nabla V(y),\mathcal
Y(y)\rangle=0$, where $\nabla V(y)$ denotes the gradient of
$V(y)$. This shows that $V(y)$ is a first integral (analytically
or formally) of system \eqref{3.1}.

Write
\[
V(y)=\sum\limits_{k=l}\limits^\infty V_k(y),\qquad \mathcal Y(y)=Ay+\sum\limits_{j=2}\limits^{\infty} G_j(y)
\]
with $l\ge 1$ a suitable natural number and $V_k(y)$ homogeneous polynomial in $y$ of degree $k$ for $k=l,l+1,\ldots$, and $G_j(y)$
vector--valued homogeneous polynomial in $y$ of degree $j$ for $j=2,3,\ldots$. Then we get from  $\langle \nabla V(y),\mathcal
Y(y)\rangle=0$ that
\begin{eqnarray}
\langle\nabla V_l,\,Ay\rangle&=&0,\label{3.2}\\
\langle\nabla
V_m,Ay\rangle&=&-\sum\limits_{j=2}\limits^m\langle\nabla
V_{m+1-j},G_j\rangle, \quad m=l+1,l+2,\ldots\label{3.3}
\end{eqnarray}
From the Bibikov's result \cite{Bi}, the linear operator $\mathcal
L_r$ from $H_r^n(y)$ to itself defined by
\[
\mathcal L_rh=\langle \nabla h(y), Ay\rangle,\qquad h(y)\in
H_r^n(y)
\]
has the spectrum $\sigma(\mathcal L_r)=\{\langle
\kappa,\lambda\rangle;\,\, \kappa\in\mathbb Z_+^n,
\,|\kappa|=r\}$. So, the solution $V_l$ of equation \eqref{3.2}
should consist of resonant monomials of degree $l$. Otherwise, it
is null by the spectrum of $\mathcal L_l$.

For each $m\in \{l+1,l+2,\ldots\}$, the right hand side of
\eqref{3.3} is an inductively known resonant polynomials of degree
$m$, because $G_j$ and $V_{m+1-j}$ are resonant homogeneous
polynomials in a vector field and in a function, respectively.
Hence it follows from the spectrum of the linear operator
$\mathcal L_m$ that $V_m$ is a resonant homogeneous polynomial of
degree $m$. This proves the claim.

By the assumption of the theorem, system \eqref{e1} has $n-1$
functionally independent analytic first integrals, denoted by
$H_1(x),\ldots,H_{n-1}$. From the Ziglin's lemma \cite{Zi} (see
also the appendix of \cite{It}), considering polynomials of these
$n-1$ first integrals with complex coefficients, we may assume
without loss of generality that the lowest order homogeneous
polynomials of these first integrals are functionally independent.

Set $V_i(y)=H_i\circ\Phi(y)$ for $i=1,\ldots,n-1$. The last claim
shows that $V_i(y)$, $i=1,\ldots,n-1$, are functionally
independent first integrals of the distinguished normal form
vector field $\mathcal Y$ of $\mathcal X$. And each $V_i(y)$
contains only resonant terms. Moreover, the lowest order
homogeneous polynomials, saying $V_i^0(y)$, of $V_i(y)$ for
$i=1,\ldots,n-1$ are also functionally independent.

The first integrals $V_i(y)$ of $\mathcal Y(y)$ satisfy the
equations $\langle \nabla V_i(y),\,Ay+g(y)\rangle=0$. Equating the
lowest order terms of these last equations, we get that $\langle
\nabla V_i^0(y),Ay\rangle=0$. Since $V_i^0(y)$ are composed of
resonant monomials, it follows that $\langle \nabla V_i^0,\lambda
y\rangle=0$ for $i=1,\ldots,n-1$, where $\lambda
y=(\lambda_1y_1,\ldots,\lambda_{n-1}y_{n-1})$. This proves that
both $Ay$ and $\lambda y$ are orthogonal to the $n-1$ dimensional
linear space spanned by $\nabla V_1,\ldots,\nabla V_{n-1}$ in an open and dense subset of $(\mathbb C^n,0)$. Hence
the vectors $Ay$ and $\lambda y$ must be parallel in an open and dense subset of $(\mathbb C^n,0)$ because we are
in the $n$--dimensional space. This implies that $Ay=\lambda y$ holds in an open and dense subset of $(\mathbb C^n,0)$ and consequently $Ay=\lambda y$ hold in $(\mathbb C^n,0)$ because they are linear in $y$.
This proves that if system \eqref{e1} has $n-1$ functionally
independent analytic first integrals, then the linear part $A$ of
\eqref{e1} should be diagonalizable.

Next we will prove that the distinguished normal form of
\eqref{e1} has the form \eqref{e2}. From the above proof we can
assume that system \eqref{e1} has the distinguished normal form of
the form $\mathcal Y
=(\lambda_1y_1+g_1(y),\ldots,\lambda_ny_n+g_n(y))$. From the
assumption of the theorem and the above proof, we know that the
vector field $\mathcal Y$ has $n-1$ functionally independent first
integrals, and each one consists of resonant polynomials. So
working in a similar way to the proof of that $A$ is diagonal, we
can verify that the two vector fields $\mathcal Y$ and $\lambda y$
should be parallel at each point $y$ in a neighborhood of the
origin. This implies that there exists a function of the form
$1+g(y)$ such that $\mathcal
Y=(\lambda_1y_1(1+g(y)),\ldots,\lambda_ny_n(1+g(y)))$.

The remainder is to prove that the distinguished normalization
from the vector fields $\mathcal X$ to $\mathcal Y$ is analytic.
For this aim we first show that if system \eqref{e1} has $n-1$
functionally independent analytic first integrals, then there
exists a $\kappa>0$ such that $|\langle
m,\lambda\rangle-\lambda_i|>\kappa$ for  all $m\in\{m\in \mathbb
Z_+^n;\,\,\langle m,\lambda\rangle-\lambda_i\ne 0, |m|\ge 2\}$.

Indeed, by Theorem 1.1 of \cite{CYZ}, i.e. the number of analytic
first integrals of system \eqref{e1} is less than or equal to
$R_\lambda$, we get that $R_\lambda= n-1$. So there exist $n-1$
linearly independent vectors $k_i=(k_{i1},\ldots,k_{in})\in\mathbb
Z_+^n$ with $|k_i|\ge 2$, $i=1,\ldots,n-1$ such that
\begin{equation}\label{3.4}
\langle k_i,\lambda\rangle=0,\quad i=1,\ldots,n-1.
\end{equation}
Since $k_1,\ldots,k_{n-1}$ are linearly independent, we can assume
without loss of generality that
\[
\det\left(\begin{array}{ccc}
k_{1,1} & \ldots & k_{1,n-1}\\
\vdots &  & \vdots\\
k_{n-1,1}& \ldots & k_{n-1,n-1}
\end{array}
\right)\ne 0.
\]
Solving \eqref{3.4} gives
\begin{equation}\label{3.5}
\lambda_1=\frac{\nu_1}{\mu_1}\lambda_n,\ldots,
\lambda_{n-1}=\frac{\nu_{n-1}}{\mu_{n-1}}\lambda_n,
\end{equation}
with $\mu_i\in \mathbb Z\setminus\{0\}$, $\nu_i\in \mathbb Z$, and
$\mu_i,\nu_i$ relatively prime for $i=1,\ldots,n-1$. We note that
$\mu_i$ and $\nu_i$ are uniquely determined by the $k_j$ for
$j=1,\ldots,n-1$. Since $\lambda\ne 0$, it follows that
$\lambda_n\ne 0$. If $\langle m,\lambda\rangle-\lambda_i\ne 0$ for
$m\in \mathbb Z_+^n$ and $|m|\ge 2$, it follows from \eqref{3.5}
that $ \left|\langle m,\lambda\rangle-\lambda_i\right|\ge
|\lambda_n|/(\mu_1\ldots\mu_{n-1})=\kappa$. This proves the claim.

This last proof shows that there does not appear the so called
small divisors in the distinguished normalization from an
analytic integrable system in $(\mathbb C^n,0)$ to its normal
form. Then working in the same way as in the proof of Lemma 2.6 of
\cite{XZ08}, we can prove that the distinguished normalization
from system \eqref{e1} to its normal form \eqref{e2} is uniformly
convergent in a neighborhood of the origin. The details are
omitted. We should mention that in our theorem part of the
eigenvalues can be zero, so for getting the coefficients of $g(y)$
in the normal form vector field $\mathcal Y$ from (2.6) of
\cite{XZ08} we must choose those $s\in\{1,\ldots,n\}$ for which
$\lambda_s\ne 0$. We complete the proof of the theorem.

\section{Proof of Theorem \ref{t3}}

\setcounter{section}{4}
\setcounter{equation}{0}\setcounter{theorem}{0}

\noindent The proof of the theorem is an improvement of that given
in \cite{XZ08} for the proof of Theorem C. Let $F$ be the
analytic integrable diffeomorphism defined on the
$n$--dimensional analytic manifold $\mathcal M$. By the
assumption, the diffeomorphism $F$ has $n-1$ functionally
independent analytic first integrals, denoted by
$V_1,\ldots,V_{n-1}$.

Let $\{U_{\alpha}\}$ be coordinate charts of $\mathcal M$ with
$\bigcup U_{\alpha}=\mathcal M$, and $x$ be the coordinate on
$U_\alpha$. Then each level surface $V_i(x)=c_i$ is invariant
under the action of $F(x)$ because by definition we have
$V_i(F(x))=V_i(x)$ for all $x\in U_\alpha$. This indicates that
each orbit of $F(x)$  is contained in
$\bigcap\limits_{i=1}\limits^n\{x\in \mathcal M;\,
V_i(x)=c_i\}:=\gamma_c$ for some $c=(c_1,\ldots,c_n)\in\mathbb
F^n$, with either $\mathbb F=\mathbb C$ or $\mathbb F=\mathbb R$.

For each $y\in U_\alpha\subset \mathcal M$, since $F$ is a
diffeomorphism on $\mathcal M$, there exists some $x\in
U_{\beta}\subset\mathcal M$ ($\beta$ may be $\alpha$ or not) such
that $y=F(x)$. Define a vector field on $\mathcal M$ by
\[
\mathcal X(y)=\det(DF(x))(\nabla V_1(F(x))\times\ldots\times\nabla
V_{n-1}(F(x))),\quad \mbox{ for } y\in\mathcal U_{\alpha},
\]
where $D$ denotes the Jacobian matrix of $F$ with respect to $x$,
$\nabla V_i(F(x))=\left.\nabla V_i(y)\right|_{y=F(x)}$ for
$i=1,\ldots,n-1$, and $\times$ denotes the cross product of
vectors in $\mathbb F^n$. We mention that the cross product of
$n-1$ vectors in $\mathbb F^n$ is defined in the proof of Lemma
\ref{l24}, and that for $v_1,\ldots,v_{n-1}\in\mathbb F^n$ their
cross product $v=v_1\times\ldots\times v_{n-1}$ is orthogonal to
each $v_i$ for $i=1,\ldots,n-1$. More generally, for $y=F^k(x)$
with some $x\in\mathcal M$ we have
\[
\mathcal X(y)=\det((DF^k)(x))(\nabla
V_1(F^k(x))\times\ldots\times\nabla V_{n-1}(F^k(x))),
\]
where $(DF^k)(x)=(DF)(F^{k-1}(x))(DF^{k-1})(x)$.

By the very definition of $\gamma_c$ and of $\mathcal X(y)$, it
follows that $\mathcal X(y)$ is an analytic vector field and is
tangent to each $\gamma_c$ at $y\in \gamma_c$. So in order for
proving $\mathcal X$ to be an embedding vector field of $F(y)$, we
only need to prove $DF(y) \mathcal X(y)=\mathcal X\circ F(y)$ for
all $y\in \mathcal M$. Because for the flow $\phi_t(y)$ of
$\mathcal X(y)$ we have $D\phi_t(y)\mathcal X(y)=\mathcal X\circ
\phi_t(y)$.

For any $y=F(x)\in\mathcal M$, since $V_i(F(y))=V_i(y)$ for
$i=1,\ldots, n-1$, we have
\begin{equation}\label{4.1}
\nabla V_i(F(y))DF(y)=\nabla V_i(y).
\end{equation}
It follows from the definition of $\mathcal X(y)$ and \eqref{4.1}
that
\begin{eqnarray}\label{4.2}
DF(y)\mathcal X(y)& =&\det(DF(x))DF(F(x))\\
&& \left(\nabla V_1(F^2(x))(DF)(F(x))\times\ldots\times \nabla
V_{n-1}(F^2(x))(DF)(F(x))\right). \nonumber
\end{eqnarray}
In addition, for $z=F(y)=F^2(x)$ and any vector $w(z)\in
T_z\mathcal M$ the tangent space of $\mathcal M$ at $z$, it
follows from the definition of cross product that
\begin{eqnarray*}
&&\left\langle w(F^2(x)),
DF(F(x))\left(\nabla V_1(F^2(x))DF(F(x))\times\ldots\times \nabla V_{n-1}(F^2(x))DF(F(x))\right)\right\rangle\\
&&\quad =\left\langle w(F^2(x))DF(F(x)),
\nabla V_1(F^2(x))DF(F(x))\times\ldots\times \nabla V_{n-1}(F^2(x))DF(F(x))\right\rangle\\
&&\quad= \det\left(\begin{array}{c} w(F^2(x))DF(F(x))\\ \nabla V_1(F^2(x))DF(F(x))\\
\vdots\\\nabla V_{n-1}(F^2(x))DF(F(x))\end{array}\right)
=\det(DF(F(x)))\det\left(\begin{array}{c} w(F^2(x))\\ \nabla V_1(F^2(x))\\
\vdots\\ \nabla V_{n-1}(F^2(x))\end{array}\right).
\end{eqnarray*}
This shows that
\begin{eqnarray}\label{4.3}
&& DF(F(x))\left(\nabla V_1(F^2(x))DF(F(x))\times\ldots\times
\nabla V_{n-1}(F^2(x))DF(F(x))\right)\nonumber\\
&&\quad =\det(DF(F(x)))\left(\nabla V_1(F^2(x))\times \ldots\times
\nabla V_{n-1}(F^2(x))\right).
\end{eqnarray}
Combining \eqref{4.2} and \eqref{4.3} we get that for $y=F(x)$
\begin{eqnarray*}
DF(y)\mathcal X(y)&=&\det(DF(x))\det(DF(F(x)))\left(\nabla
V_1(F^2(x))\times \ldots\times \nabla
V_{n-1}(F^2(x))\right)\\
&=&\det (DF^2(x))\left(\nabla V_1(F^2(x))\times \ldots\times
\nabla
V_{n-1}(F^2(x))\right)\\
&=&\mathcal X\circ F(y).
\end{eqnarray*}
This shows that $\mathcal X(y)$ is an embedding vector field of
$F(y)$. We complete the proof of the theorem.

\end{document}